\documentstyle{amltd}
\begin{document}
 
\annalsline{158}{2003}
\received{February 23, 2001}
\startingpage{165}
\def\bye{\end{document}}
 \font\tenrm=cmr10
\def\ritem#1{\item[{\rm #1}]}
\catcode`\@=11
\font\twelvemsb=msbm10 scaled 1100
\font\tenmsb=msbm10
\font\ninemsb=msbm10 scaled 800
\newfam\msbfam
\textfont\msbfam=\twelvemsb  \scriptfont\msbfam=\ninemsb
  \scriptscriptfont\msbfam=\ninemsb
\def\msb@{\hexnumber@\msbfam}
\def\Bbb{\relax\ifmmode\let\next\Bbb@\else
 \def\next{\errmessage{Use \string\Bbb\space only in math
mode}}\fi\next}
\def\Bbb@#1{{\Bbb@@{#1}}}
\def\Bbb@@#1{\fam\msbfam#1}
\catcode`\@=12

 \catcode`\@=11
\font\twelveeuf=eufm10 scaled 1100
\font\teneuf=eufm10
\font\nineeuf=eufm7 scaled 1100
\newfam\euffam
\textfont\euffam=\twelveeuf  \scriptfont\euffam=\teneuf
  \scriptscriptfont\euffam=\nineeuf
\def\euf@{\hexnumber@\euffam}
\def\frak{\relax\ifmmode\let\next\frak@\else
 \def\next{\errmessage{Use \string\frak\space only in math
mode}}\fi\next}
\def\frak@#1{{\frak@@{#1}}}
\def\frak@@#1{\fam\euffam#1}
\catcode`\@=12

\def\ra{\rightarrow} \def\Ra{\Rightarrow}
\def\rhu{\rightharpoonup} \def\rhd{\rightharpoondown}
\def\asmor{\hbox{$\rightarrow\hskip-7pt\rightarrow$}}
\def\dst{\displaystyle} 
\def\ie{{i.e.\ }} \def\eg{{e.g.\/}\ }
\def\etc{{\it etc.\/}\ }
\def\pt{{\rm point}}\def\ev{{\rm ev}}
\def\C{\bf C}
\def\Z{\bf Z}
\def\D{{\cal D}}

\def\og{\leavevmode\raise
.3ex\hbox{$\scriptscriptstyle\langle\!\langle$}}
\def\fg{\leavevmode\raise
.3ex\hbox{$\scriptscriptstyle\,\rangle\!\rangle$}}

\newcount\notenumber \notenumber=1
\def\note#1{\footnote{($^{\the\notenumber}$)}{#1}\global%
\advance\notenumber by 1}

\font\pcap=cmcsc10 \font\bfit=cmbxti10
\font\twelvebf=cmbx12
\font\sevenrm=cmr7
\font\eightrm=cmr8
 
\def\hb{\hfill\break}
 
\def \H{{\cal H}}
\def \E{{\cal E}}
\def \F{{\cal F}}
\def \A{{\cal A}}
\def \B{{\cal B}}
\def \cS{{\cal S}}
\def \L{{\cal L}}
\def \K{{\cal K}}
\def \M{{\cal M}}
\def \R{{\bf R}}
\def \N{{\bf N}}
\def \C{{\bf C}}
\def \Cst{{$C^*$}}
\def\sd{\hbox{$\times\!$\vbox{\hrule height 5pt width .5pt}}\,}
\def \mt{\mapsto}

\font\tengoth=eufm10
\font\sevengoth=eufm7
\font\fivegoth=eufm5
\newfam\gothfam
\textfont\gothfam=\tengoth
\scriptfont\gothfam=\sevengoth
\scriptscriptfont\gothfam=\fivegoth
\def\goth{\fam\gothfam\tengoth}
\def\gtM{{\frak M}}

\def\pt{{\rm point}}

\def\bcl{$\widehat {\rm BC}^{\hbox{\sevenrm lim}}_{\hbox{\sevenrm
coef}}$}
\def\bc{$\widehat {\rm BC}_{\hbox{\sevenrm coef}}$}

\def\bcl{$\widehat {\rm SNC}^{\hbox{\sevenrm lim}}_{\hbox{\sevenrm
coef}}$}
\def\bc{$\widehat {\rm SNC}_{\hbox{\sevenrm coef}}$}

\def\lra{\longrightarrow}

\def \limind{\mathop{\oalign{lim\cr
\hidewidth$\longrightarrow$\hidewidth\cr}}}
\def\beg{\underline{E}\Gamma }

\title{Groups acting properly on ``bolic'' spaces\\ and the Novikov conjecture} 
\shorttitle{``Bolic" spaces and the Novikov conjecture}  
  \twoauthors{Gennadi Kasparov}{Georges Skandalis}
 \institutions{Institut de Math\'ematiques de Luminy,  Marseille, France\\
{\eightpoint {\it E-mail address\/}: kasparov@iml.univ-mrs.fr}\\
\vglue6pt
Institut de Math\'ematiques de Jussieu, Universit\'e Denis Diderot 
(Paris VII), Paris, France\\
{\eightpoint {\it E-mail address\/}: skandal@math.jussieu.fr
}}

\centerline{\bf Abstract}
\vglue12pt
We introduce a class of metric spaces which we call ``bolic''. They include
hyperbolic spaces, simply connected complete manifolds of nonpositive
curvature, euclidean buildings, etc. We prove the Novikov conjecture on
higher signatures for any discrete group which admits a proper isometric
action on a ``bolic'', weakly geodesic metric space of bounded geometry.

\section{Introduction}

This work has grown out of an attempt to give a purely $KK$-theoretic
proof of a result of A.\ Connes and H.\ Moscovici ([CM], [CGM]) that
hyperbolic groups satisfy the Novikov conjecture. However, the main
result of the present paper appears to be much more general than this.
In the process of this work we have found a class of metric spaces
which contains hyperbolic spaces (in the sense of M.\ Gromov), simply
connected complete Riemannian manifolds of nonpositive sectional
curvature, euclidean buildings, and probably a number of other
interesting geometric objects. We called these spaces ``bolic spaces''.

Our main result is the following:

\proclaim{Theorem} Novikov\/{\rm '}\/s conjecture on {\rm ``}\/higher signatures\/{\rm ''} is true for any discrete
group acting properly by isometries on a weakly bolic{\rm ,} weakly geodesic
metric space of bounded coarse geometry.
\endproclaim
\begin{itemize}
\ritem{--} The notion of a ``bolic'' and ``weakly bolic'' space is defined in
Section 2, as well as the notion of a ``weakly geodesic'' space;
 
\ritem{--} bounded coarse geometry (\ie bounded
geometry in the sense of P. Fan; see [HR]) is discussed in Section 3.
\end{itemize}
All conditions of the theorem are satisfied, for example, for any
discrete group acting properly and isometrically either on a simply
connected complete\break Riemannian manifold of nonpositive, bounded sectional
curvature, or on a euclidean building  \pagebreak with uniformly bounded
ramification numbers. All conditions of the theorem are also satisfied
for word hyperbolic groups, as well as for finite products of groups of
the above classes. Note also that the class of
geodesic bolic metric spaces of bounded geometry is closed under taking
finite products (which is not true, for example, for the class of
hyperbolic metric spaces). 

The Novikov conjecture for discrete groups which belong to the above
described classes was already proved earlier by different methods. In the
present paper we give a proof valid for all these cases simultaneously,
without any special arrangement needed in each case separately.
Moreover, the class of bolic spaces is not a union of the above classes
but probably is much wider. Although we do not have at the moment any
new examples of bolic spaces interesting from the point of view of the
Novikov conjecture, we believe they may be found in the near future.

 In [KS2] we announced a proof of the Novikov
conjecture for discrete groups acting properly, by isometries
on geodesic {\it uniformly locally finite}\/ bolic metric spaces. The
complete proof was given in a preprint, which remained unpublished
since we hoped to improve the uniform local finiteness condition. This is
done in the present paper where uniform local finiteness is replaced by
a much weaker condition of bounded geometry.

 Our proof follows the main lines of [K2] and [KS1]: we construct a
`proper' $\Gamma $-algebra $\A$, a `dual Dirac' element $\eta\in
KK^\Gamma (\C,\A)$ and a `Dirac' element in $KK^\Gamma (\A,\C)$. In
the same way as in [KS1], the construction of the dual Dirac element relies
on the construction of an element $\gamma\in KK_\Gamma (\C,\C)$ (the
Julg-Valette element in the case of buildings; cf.\ [JV]).

 Here is an explanation of the construction of
these ingredients:

The algebra $\A$ is constructed in the following way (\S7): We may
assume that our bolic metric space $X$ is locally finite (up to
replacing it by a subspace consisting of the preimages in $X$ of the
centers of balls of radius $\delta$ covering $X/\Gamma $). The
assumption of bounded geometry is used to construct a `good'\break  $\Gamma
$-invariant measure $\mu$ on
$X$. Corresponding to the Hilbert space $\H=L^2(X,\mu)$ is a \Cst-algebra
$\A(\H)$ constructed in [HKT] and [HK]; denote by $H$ the subspace of
$\Lambda ^*(\ell^2(X))$ spanned by $e_{x_1}\wedge\cdots\wedge e_{x_p}$,
where the set $\{x_1,\ldots ,x_p\}$ has diameter $\le N$ (here $N$ is a
large constant appearing in our construction and related to bolicity);
then $\A$ is a suitable proper subalgebra of $\K(H) \,\widehat{\otimes}\,\A(\H)$.

 The inclusion of $\A$ in $\K(H) \,\widehat{\otimes}\,\A(\H)$ together with the
Dirac element of $\A(\H)$ constructed in [HK], gives us the Dirac element
for $\A$.

 The element $\gamma$ (\S6) is given by an operator $F_x$ acting on
the Hilbert space $H$ mentioned above, where $x\in X$ is a point chosen
as the origin. The operator $F_x$ acts on $e_{x_1}\wedge\cdots\wedge
e_{x_p}$ as Clifford multiplication by a unit vector $\phi_{S,x}\in
\ell^2(X)$ where
$S=\{x_1,\ldots ,x_p\}$ and $\phi_{S,x}$ has support in a set $Y_{S,x}$ of
points closest to $x$ among the points in $S$ or
points which can be added to $S$ keeping the diameter of $S$ not
greater than $N$. The bolicity condition is used here. Namely:

We prove that if $y\in Y_{S,x}$, denoting by $T$ the
symmetric difference of $S$ and $\{y\}$, we have $\phi_{S,x}= \phi
_{T,x}$, which gives that $F_x^2-1\in \K(H)$ (this uses half of
the bolicity, namely condition (B2$'$)).

Averaging over the radius of a ball centered at $x$   used
in the construction of $\phi_{S,x}$ allows us to prove that $\lim_{S\ra \infty}
\|\phi_{S,x}-\phi_{S,y}\|=0$, whence $F_x-F_y\in \K(H)$ for any $x,y\in
X$, which shows
that $F_x$ is $\Gamma $-invariant up to $\K(H)$ (this uses condition (B1)).

  In the same way as $\phi_{S,x}$, we construct a measure
$\theta_{S,x}$ supported by the points of $S$ which are the most remote
from $x$. This is used as the center for the `Bott element' in the
construction of the dual Dirac element (Theorem 7.3.a).  

There are also some additional difficulties we have to deal with:

\begin{itemize}
\ritem{a)} Unlike the case of buildings (and the hyperbolic case), we do not
know anything about contractibility of the Rips complex. We need to use an
inductive limit argument, discussed in Sections 4 and 5.

\ritem{b)} The Dirac element appears more naturally as an
asymptotic $\Gamma $-morphism. On the other hand, since we wish to obtain
the injectivity of the Baum-Connes map in the reduced $C^*$-algebra,
we need to use $KK$-theory. This is taken care of in Section 8.
\end{itemize}

Our main result on the Novikov conjecture naturally corresponds to the injectivity part 
of the Baum-Connes conjecture for the class of groups that we consider
(see Theorem 5.2). 
We do not discuss the surjectivity part of the Baum-Connes conjecture
(except maybe in Proposition 5.11). We can mention however that our result has already been used by 
V. Lafforgue in order to establish the Baum-Connes conjecture for a certain 
class of groups ([L]). On the other hand, M.\ Gromov has recently given ideas for construction 
of examples of discrete groups which do not admit any uniform embedding into a Hilbert space ([G1], [G2]). For these groups the surjectivity part of the Baum-Connes conjecture 
{\it with coefficients}\/ fails ([HLS]).

 The paper is organized as follows: in Sections 2--4 we introduce
the main definitions. Sections 5--8 contain the mains steps of the proof.
More precisely: \begin{itemize}
\ritem{--} Bolicity is defined in Section 2, where we prove that hyperbolic spaces
and Riemannian manifolds of nonpositive sectional curvature are bolic.
 
\ritem{--} The property of bounded geometry is discussed in Section 3.
 
\ritem{--} Section 4 contains some
preliminaries on universal proper $\Gamma $-spaces and Rips complexes.
 
\ritem{--} Section 5 gives the statement of our main result and a general
framework of the proof.
 
\ritem{--} Section 6 contains the construction of the $\gamma$-element.
 
\ritem{--} Finally, in Sections 7 and 8 we explain the construction of the
\Cst-algebra of a Rips complex, give the construction of the dual Dirac
and Dirac elements in $KK$-theory, and finish the proof of our main result.
\end{itemize}

The reader is referred to [K2] for the main definitions related to the equivariant $KK$-theory,  graded
algebras, graded tensor products and for some related jargon: for example, $\Gamma $-algebras are just
\Cst-algebras equipped with a continuous action of a locally compact group
$\Gamma $, $C(X)$-algebras are defined in [K2], 1.5, etc. Unless otherwise
specified, all tensor products of
\Cst-algebras are considered with the minimal \Cst-norm.
{\it All groups acting on \Cst{\rm -}algebras are supposed to be locally compact and
$\sigma$\/{\rm -}\/compact{\rm ,} all discrete groups {\rm --} countable.}

\section{``Bolicity''} 

Let $\delta$ be a nonnegative real number. Recall that a map (not
necessarily continuous) $f:X \ra X'$ between metric spaces
$(X,d)$ and $(X',d')$ is said to be a {\it $\delta$-isometry}\/
if for every
pair $(x,y)$ of elements of $X$ we have $ |d'(f(x),f(y)) - d(x,y)|
\le \delta$. Also, the metric space $(X,d)$ is said to be
$\delta$-geodesic if for every pair $(x,y)$ of points of $X$, there
exists a $\delta$-isometry $f:[0,d(x,y)] \ra X$ such that $f(0)=
x\,,\;f(d(x,y)) = y$.

\numbereddemo{Definition} The space $(X,d)$ is said to be {\it weakly
$\delta$\/{\rm -}\/geodesic}\/ if for every pair $(x,y)$ of points of $X$, and
every $t\in[0,d(x,y)]$ there exists a point $a\in X$ such that
$d(a,x)\le t+\delta$ and $d(a,y)\le d(x,y)-t+\delta$. The point
$a\in X$ is said to be a {\it $\delta$\/{\rm -}\/middle point}\/ of $x,y$ if
$|2d(x,a)-d(x,y)| \le 2\delta$ and $|2d(y,a)-d(x,y)| \le 2\delta$.
We will say that the space $(X,d)$ admits $\delta$-middle points if there
exists a map $m:X\times X\ra X$ such that for any $x,y\in X$, the point
$m(x,y)$ is a $\delta$-\/{\it middle point} of $x,y$. The map $m$ will be called a
$\delta$-middle point map.
\enddemo

Note that in the above definition of a weakly $\delta$-geodesic space, one
can obviously take $t\in [-\delta,0]\cup [d(x,y) ,d(x,y)+\delta]$ and
$a=x$ or $a=y$. This will be useful in Section 6.
Also note that a $\delta$-geodesic space is weakly
$\delta$-geodesic. In a weakly
$\delta$-geodesic space, every pair of points admits a $\delta$-middle
point.

\numbereddemo{Definition} We will say that a metric space $(X,d)$ is
{\it $\delta $\/{\rm -}\/bolic}\/ if:
\begin{itemize}
\item[(B1)] For all 
$r>0$, there exists 
$R>0$ such that for every quadruple $x,y,z,t$ of points of $X$
satisfying
$d(x,y)+d(z,t)\le r$ and $d(x,z)+d(y,t) \ge R$, we have
$d(x,t)+d(y,z) \le d(x,z)+d(y,t)+2\delta $. 
\item[(B2)]
There exists a map $m:X\times X\ra X$ such that for all $x,y,z\in
X$ we have $ 2d(m(x,y),z) \le \left
(2d(x,z)^2+2d(y,z)^2-d(x,y)^2\right)^{1/2}+4\delta$.
\end{itemize}
\enddemo

 We will say that a metric space
$(X,d)$ is {\it weakly $\delta $\/{\rm -}\/bolic} if it satisfies the condition
(B1) and the following condition:
\begin{itemize}
\item[(B2$'$)] There exists a $\delta$-middle point map
$m:X\times X\ra X$ such that if $x,y,z$ are points of $X\,,$ then
$d(m(x,y),z)<\max
(d(x,z),d(y,z)) + 2\delta$. Moreover, for every $p\in {\bf R}_+$,
there exists $N(p)\in {\bf R}_+$ such that for all $N\in {\bf
R}_+$,\break $N\ge N(p)$, if
$d(x,z)\le N$, $d(y,z)\le N$ and
$d(x,y)> N$ then $d(m(x,y),z)\break<N-p$.
\end{itemize}

Condition (B2$'$) is a property of ``strict convexity'' of balls. Bolic
spaces are obviously weakly bolic (a point $m(x,y)$ satisfying
condition (B2) is automatically a $2\delta$-middle point of $x,y$;
apply condition (B2) to $z=x$ and
$z=y$).

\proclaim{Proposition} Any $\delta$-hyperbolic space admitting $\delta$\/{\rm -}\/middle points is
$3\delta/2$\/{\rm -}\/bolic.
\endproclaim

\demo{Proof} Let $(X,d)$ be a $\delta $-hyperbolic metric space. Condition
(B1) is obviously satisfied.

Assume moreover that we have a $\delta$-middle point map
$m:X\times X\ra X$. Let $z\in X$. The hyperbolicity condition gives:
\begin{eqnarray*}&&
 d(z,m(x,y))+d(x,y)\\[4pt]
 &&\qquad \le
\sup\,\{\,d(y,z)+d(x,m(x,y))\,,\; d(x,z)+d(y,m(x,y))\,\}+2\delta
\\[4pt]
&&\qquad\le \sup\,\{\,d(x,z)\,,\, d(y,z)\,\}+{\dst d(x,y)+2\delta\over
2}+2\delta\,.
\end{eqnarray*}
 Therefore, $$2d(z,m(x,y)) \le\, 2\sup\,\{\,d(x,z)\,,\; d(y,z)\,\}-
d(x,y) +6\delta\,.$$ Now, if $s,t,u$ are nonnegative real numbers such
that
$\vert t-u\vert \le s$, we have $$(2t-u)^2+u^2=2t^2+2(t-u)^2\le
2t^2+2s^2\,.$$ Setting $s=\inf\,\{\,d(x,z)\,,\;
d(y,z)\,\}\,,\;t=\sup\,\{\,d(x,z)\,,\; d(y,z)\,\}$ and $u=d(x,y)$, we
find
\vglue12pt
\noindent \hfill ${\displaystyle 2\sup\,\{\,d(x,z)\,,\; d(y,z)\,\}- d(x,y)\le \left
(2d(x,z)^2+2d(y,z)^2-d(x,y)^2\right)^{1/2}\,.}$
\enddemo

\proclaim{Proposition} Every nonpositively curved simply connected complete Riemannian
manifold is $\delta$\/{\rm -}\/bolic for any $\delta>0$.
\endproclaim

In particular Euclidean spaces, as well as symmetric spaces $G/K$,
where $G$ is a semisimple Lie group and $K$ its maximal
compact subgroup, are bolic.

\demo{Proof} Let us first prove (B2). Recall the cosine theorem for nonpositively
curved manifolds (cf.\ [H, 1.13.2]): For any geodesic triangle with edges
of length $a,\;b$ and $c$ and the angle between the edges of the length
$a$ and
$b$ equal to $\alpha$, one has:
$$a^2+b^2-2ab\cos\alpha\le c^2.$$ Define $m(x,y)$ as the middle point of
the unique geodesic segment joining $x$ and $y$. Apply the cosine
theorem to the two geodesic triangles: $(x,z,m(x,y))$ and
$(y,z,m(x,y))$. If we put
$a=d(x,z),\;b=d(y,z),\;c=d(x,m(x,y))=d(y,m(x,y)),\;e=d(z,m(x,y))$ then
$$c^2+e^2-2ce\cos\alpha\le a^2,\;\;c^2+e^2-2ce\cos(\pi -\alpha )\le
b^2$$ where the angle of the first triangle opposite to the edge $(x,z)$
is equal to
$\alpha$. The sum of these two inequalities gives (B2) with
$\delta =0$.

For the proof of (B1), let $x$ and $y\in X$. Suppose that $z(s),\;0\le
s\le d(z,t),$ is a geodesic segment (parametrized by distance) joining
$t=z(0)$ with $z=z(d(z,t))$. Then it follows from the cosine
theorem that
$$\vert (\partial /\partial s)(d(y,z(s))-d(x,z(s)))\vert\le {\dst 2c\over
{a(s)+b(s)}},$$ where $c=d(x,y),\;a(s)=d(x,z(s)),\;b(s)=d(y,z(s))$.

Indeed, the norm of the derivative on the left-hand side does not exceed
$\|{\rm grad} f(u)\|$, where $f(u)=d(x,u)-d(y,u)$ is a function of $u=z(s)$.
It is clear that $\|{\rm grad} f(u)\|$ is the norm of the difference between
the two unit vectors tangent to the geodesic segments $[x,u]$ and
$[y,u]$ at the point $u$, so that  $\|{\rm grad} f(u)\|^2=2(1-\cos \alpha)$, where
$\alpha$ is the angle between these two vectors. The cosine theorem applied
to the geodesic triangle $(x,y,u=z(s))$ gives: $a(s)^2+b(s)^2-c^2\le
2a(s)b(s)\cos
\alpha$, whence $2a(s)b(s)(1-\cos
\alpha)\le c^2-(a(s)-b(s))^2$. Therefore, $$\|{\rm grad} f(u)\|^2\le
{c^2-\big(a(s)-b(s)\big)^2\over a(s)b(s)}\le {4c^2 \over \big(a(s)+b(s)\big)^2}$$
since $c\le
a(s)+b(s)$. This implies the above inequality.

Integrating this inequality over
$s$, one gets the estimate:
$$(d(y,z)-d(x,z))-(d(y,t)-d(x,t))\le {\dst 2\over R-r}d(x,y)d(z,t) \eqno
(1)$$ with $R$ and $r$ as in the condition (B1), which gives (B1) with
$\delta$ arbitrarily small.
\enddemo

\proclaim{Proposition} Euclidean buildings are $\delta$\/{\rm -}\/bolic for any $\delta>0$.
\endproclaim

\demo{Proof} The property (B2) (with $\delta=0$) is proved in [BT, Lemma 3.2.1].
To prove (B1) let us denote the left side of (1) by $q(x,y;z,t)$. Then,
clearly,
$q(x,y;z,t)+q(y,u;z,t)=q(x,u;z,t)$. The same type of additivity holds also
in the $(z,t)$-variables. Now when the points $(x,y)$ are in one chamber
and points $(z,t)$ in another one, we can apply the inequality (1) because in
this case all four points $x,y,z,t$ belong to one apartment. In general we
reduce the assertion to this special case by using the above additivity
property. 
\enddemo

\proclaim{Proposition}
A product of two bolic spaces when endowed with the distance
such that $d((x,y),(x',y'))^2=d(x,x')^2+d(y,y')^2$ is bolic.
\endproclaim

\demo{Proof} Let $(X_1,d)$ and $(X_2,d)$ be two $\delta$-bolic spaces. We
show that\break $X_1\times X_2$ is $2\delta$-bolic. Take
$r>0$ and let $R$ be the corresponding constant in the condition (B1) for
both $X_i$. Let $R'\in \R_+$ be big enough. For $x_i,y_i,z_i,t_i\in X_i$, put
$x=(x_1,x_2)\,,\;y=(y_1,y_2)\,,\; z=(z_1,z_2)$ and $t=(t_1,t_2)$.
Assume that $d(x,y)+d(z,t)\le r$ and
$d(x,z)+d(y,t)\ge R'$. We distinguish two cases:
\begin{itemize}
\item[--] We have $d(x_1,z_1)+d(y_1,t_1)\ge R$ and
$d(x_2,z_2)+d(y_2,t_2)\ge R$.
\end{itemize}
In this case
$$d(x_i,t_i)+d(y_i,z_i)\le d(x_i,z_i)+d(y_i,t_i)+2\delta.$$ Put
$$z'_i=d(x_i,z_i),\;y'_i=d(x_i,t_i)-d(y_i,t_i)\,,\;t'_i=d(x_i,t_i).$$ Note
that $$d(y_i,z_i)\le
d(x_i,z_i)+d(y_i,t_i)-d(x_i,t_i)+2\delta=z'_i-y'_i+2\delta.$$ Note also
that
 $|y'_i|\le d(x_i,y_i)$ and $|z'_i-t'_i|\le d(z_i,t_i).$  Put $x'=(0,0)$,
$y'=(y'_1,y'_2)$,
$z'=(z'_1,z'_2)$ and
$t'=(t'_1,t'_2)$. As $\R^2$ is $\delta'$-bolic for every $\delta'$, if
$R'$ is large enough, we find that $$\|z'-y'\|+\|t'-x'\|\le
\|z'-x'\|+\|t'-y'\|+(4-2\sqrt2)\delta.$$ Now
$\|z'-x'\|=d(x,z)\,,\;\|t'-x'\|=d(x,t)\,,
\;\|t'-y'\|=d(y,t)$ and
$ d(y,z)\le\|y'-z'\|+2\sqrt 2\delta$. We therefore get condition (B1) in
this case.
\begin{itemize}
\item[--] We have $d(x_2,z_2)+d(y_2,t_2)\ge R$ but
$d(x_1,z_1)+d(y_1,t_1)\le R$. 
\end{itemize}
Choosing $R'$ large
enough, we may assume that if $s,u\in \R_+$ are such that $s\le R+r$ and
$(s^2+u^2)^{1/2}\ge R'/2-r$, then $(s^2+u^2)^{1/2}\le u+\delta$.
Therefore, $d(y,z)\le d(y_2,z_2)+\delta$ and $d(x,t)\le
d(x_2,t_2)+\delta$, whence condition (B1) follows also in this case.

Let us check condition (B2). Let $x_1,y_1,z_1\in X_1$ and
$x_2,y_2,z_2\in X_2$. Put $A_i=\left
(2d(x_i,z_i)^2+2d(y_i,z_i)^2-d(x_i,y_i)^2\right)^{1/2}$ ($i=1,2$). We
have 
\begin{eqnarray*}
4(d(m_1(x_1,y_1),z_1)^2+d(m_2(x_2,y_2),z_2)^2)&\le&
(A_1+4\delta)^2+(A_2+4\delta)^2\\
&\le&((A_1^2+A_2^2)^{1/2}+4\sqrt 2\delta)^2
\end{eqnarray*}
and condition (B2)
follows.
\enddemo

\numbereddemo{{R}emark} Let $X$ be a $\delta$-bolic space, and let $Y$ be a
subspace of
$X$ such that for every pair $(x,y)$ of points of $Y$ the
distance of $m(x,y)$ to $Y$ is $\le\delta$. Then $Y$ is
$2\delta$-bolic. The same is true for weakly bolic spaces.
\enddemo

{\it {R}emark} 2.8. Bolicity is very much a euclidean condition. On the other hand, weak
bolicity, is not at all euclidean. Let $E$ be a finite-dimensional normed
space.
\begin{itemize}
\item[(a)] If the unit ball of the dual space $E'$ is strictly convex then $E$
satisfies
condition (B1). 
\item[(b)] If there are no segments of length $1$ in the unit sphere of $E$, then $E$
satisfies condition (B2$'$).
\end{itemize}
 \advance\theoremcount by 1
Indeed, an equivalent condition for the strict convexity of the unit ball
of $E'$ is
that for any nonzero $x\in E$, there exists a unique $\ell_x$ in the unit
sphere of $E'$ such that $\ell_x(x)=\|x\|$; moreover, the map
$x\mapsto \|x\|$
is differentiable at $x$, its differential is $\ell_x$ and the map
$x\mapsto
\ell_x$ is continuous and homogeneous (\ie $\ell _{\lambda
x}=\ell_x$ for
$\lambda>0$).

 Now, let $r>0$. There exists an $\varepsilon>0$ such that for all
$u,v\in E$ of norm $1$, if $\|u-v\|\le \varepsilon\,,$ then
$\|\ell_u-\ell_v\|\le \delta/r$. Take $x,y,z,t\in E$
satisfying $\|x-y\|\le r\,,\;\|z-t\|\le r$ and $\|x-z\|\ge
2r/\varepsilon +r $. Note that for nonzero $u,v\in
E$, we have $\|\|u\|^{-1}u-\|v\|^{-1}v\|\le 2\|u-v\|\|u\|^{-1}$.

For every $s\in [0,1]\,,$ set $x_s=sx+(1-s)y$. Since $\|x_s-z\|\ge
2r/\varepsilon$, the distance between $u_s=\|x_s-z\|^{-1}(x_s-z)$ and
$v_s=\|x_s-t\|^{-1}(x_s-t)$ is $\le \varepsilon$. Therefore the
derivative of $s\mapsto \|x_s-z\|-\|x_s-t\|$, which is equal to
$(\ell_{u_s}- \ell_{v_s}) (x-y)$, is $\le \delta$. Therefore condition
(B1) is satisfied.

Assume now that there are no segments of length $1$ in the unit
sphere of~$E$. Let $k=\sup\{\|y+z\|/2\,,\;\|y\|\le 1\,,\;\|z\|\le
1\,\;\|y-z\|\ge 1\,\}$. By compactness and since there are no segments of
length $1$ in the unit sphere of $E$, $k<1$. If $x,y,z\in E$ satisfy
$\|x-z\|\le N\,,\;\|y-z\|\le N\,,$ and $\|x-y\|\ge N$, then
$\|z-(x+y)/2\|\le kN$. Setting $m(x,y)=(x+y)/2$ we obtain condition (B2$'$)
because for any $p>0$ there is an $N>0$ such that $kN<N-p$.  

\numbereddemo{{R}emark} It was proved recently by M. Bucher and A. Karlsson ([BK]) 
that condition (B2) actually implies (B1). 
\enddemo
\vglue-8pt

\section{Bounded geometry}
\vglue-4pt
Consider a metric space $(X,d)$ which is proper in the sense that any
closed bounded subset in $X$ is compact.
Let us fix some notation:

For $x\in X$ and $r\in {\bf R}_{+}$, let $B(x,r)=\{\,y\in
X\,,\;d(x,y)< r\,\}$ be the open ball with center $x$ and radius $
r$ and $\overline {B}(x,r)=\{\,y\in X\,,\;d(x,y)\le r\,\}$ the closed
ball
with center $x$ and radius $ r$.

The following condition of {\it bounded coarse geometry} will be
important for us. Recall from [HR] its definition:

\numbereddemo{Definition} A metric space $X$ has {\it bounded coarse geometry}\/ if there exists 
$\delta >0$ such that  for any $R>0$ there exists $K=K(R)>0$ such that in any
closed ball of radius $R$, the maximal number of points with pairwise
distances between them $\ge\delta$ does not exceed $K$.
\enddemo

We need to consider a situation in which a locally compact group $\Gamma$ acts
properly by isometries on $X$. For simplicity we will assume in this section
that $\Gamma$ is a {\it discrete} group.

\proclaim{Proposition}
Let $X$ be a proper metric space of bounded
coarse geometry and $\Gamma$ a discrete group which acts properly and
isometrically on $X$. Then there exists on $X$ a
$\Gamma $\/{\rm -}\/invariant positive measure $\mu$ with the property
that  for any
$R>0$ there exists $K>0$ such that for any $x\in
X\,,\;\mu(\overline B(x,R))\le K$ and $\mu(B(x,2\delta ))\ge 1$.
\endproclaim

\demo{Proof} Let $Y$ be a maximal subset of points of $X$ such that the
distance between any point of $Y$ and a $\Gamma $-orbit passing
through any other point of $Y$ is $\ge\delta$; by maximality of
$Y$, for any $x\in X,\;\,d(x,\,\Gamma \cdot Y)<\delta$. For $y\in
Y$, let $n(y,\delta )$ be the number of points of $\Gamma y\cap
B(y,\delta)$. Define a measure on $X$ by assigning to any point on
the orbit $\Gamma y$ the mass
$n(y,\delta )^{-1}$. In this way we define a $\Gamma$ invariant measure
$\mu$ on the set
$\Gamma \cdot Y$. Outside of this set, put $\mu$ to be
$0$. Note that for any $z\in \Gamma \cdot
Y\,,\,\;\mu(B(z,\delta))=1$.

For any $x\in X$, there exists $z\in \Gamma \cdot Y$ such that
$d(x, z)< \delta$; hence $\mu(B(x,2\delta ))\ge
\mu(B(z,\delta ))=1$.

For any $x\in X$ and $R>0$, let $Z$ be a maximal subset of
$\Gamma \cdot Y\cap \overline B(x,R)$, with pairwise distances between any
two points
$\ge\delta$. By definition, $Z$ has at most $K(R)$ points.
Obviously $\Gamma \cdot Y\cap \overline B(x,R)\subset \dst\bigcup_{z\in Z}
B(z, \delta)$; therefore 
\vglue12pt
\noindent   ${\displaystyle \mu(\overline B(x,R))=\mu(\Gamma \cdot Y\cap
\overline
B(x,R)) \le \mu(\dst\bigcup_{z\in Z} B(z,
\delta))\le \sum_{z\in Z} \mu(B(z, \delta))\le K(R)\,.}$
\enddemo
\vglue12pt
The following converse to the above proposition can be used in order to
give examples of bounded coarse geometric spaces.
 
\proclaim{Proposition}
Assume that $X$ is a metric space equipped with a positive
measure $\mu$ {\rm (}\/not necessarily $\Gamma $\/{\rm -}\/invariant\/{\rm )} which satisfies the
following condition\/{\rm :} there exists $\delta$ such that for all $R>0${\rm ,} 
there exists $\tilde {K}=\tilde {K}(R)>0$ such that for any $x\in
X\,,\;\mu(B(x,R))\le\tilde {K}$ and $\mu(B(x,\delta /2))\ge 1$. Then $ X$
is a bounded coarse geometric \pagebreak space.
\endproclaim

\demo{Proof} Let $y_1,\ldots ,y_p\in \overline B(x,R)$ be points with pairwise distances $
\ge\delta$. Then the balls $B(y_1,\delta /2),\ldots ,B(y_p,\delta /2)$ do
not intersect and are all contained in $B(x,R+\delta /2)$. Therefore,
according to our assumption, $\tilde {K}(R+\delta /2)\ge
\mu(B(x,R+\delta /2))\ge p$.
\enddemo

We will call a discrete metric space $(X,d)$ {\it locally finite} if any
ball contains only a finite number of points.

\numbereddemo{{R}emark} All locally finite metric spaces equipped with an isometric proper
action of a discrete group $\Gamma$, which have only a finite number of
orbits of $\Gamma $-action, have bounded coarse geometry (BCG). All complete
Riemannian manifolds with sectional curvature bounded from below are
BCG-spaces. (This follows from Rauch's comparison theorem together with
the criterion given in Proposition 3.3, the measure
$\mu$ is the one defined by the Riemannian metric.) Euclidean
buildings with uniformly bounded ramification numbers are\break BCG-spaces.
A finite product of BCG-spaces is a BCG-space.
Bounded coarse geometry is obviously hereditary with respect to
passing to subspaces.\break Together with the hereditary property of bolicity
(see Remark 2.7 of the previous section), this gives a large number of
examples of {\it locally finite} bolic metric spaces of bounded coarse
geometry. We record this for future use in the following:

\proclaim{Proposition}
In any bolic{\rm ,} weakly geodesic metric space of bounded coarse
geometry equipped with an isometric proper action of a discrete group
$\Gamma ${\rm ,} there exists a
$\Gamma $\/{\rm -}\/invariant{\rm ,} locally finite{\rm ,} bolic{\rm ,} weakly geodesic metric
subspace of bounded coarse geometry. The assertion remains true if we
replace bolicity by weak bolicity.
\endproclaim

\vglue-20pt
 
\section{Rips complexes}
\vglue-4pt

Before we state (in the next section) our main result, we would like to
introduce one more technical tool which will play a crucial role in the
proof. Recall from [BCH] that there exists a ``universal example'' $\beg $ for
proper actions of a locally compact group $\Gamma $. We will give now its
construction in a form suitable for our purposes.

Let $X$ be a locally compact metrizable $\sigma$-compact space. We will
denote by
$\gtM$ the set of finite positive measures on $X$ with total mass
contained in
$(1/2,1]$, endowed with the topology of duality with the algebra of
continuous functions with compact support. Clearly,
$\gtM=K-{1\over2}K$ where $K$ is the set of finite positive
measures on $X$ with total mass $\le1$. As $K$ is compact
and $\gtM$ is open in $K$, $\gtM$ is locally compact.

Let $\Gamma $ be a locally compact group acting properly on the
space $X$. Then $\Gamma $ acts naturally on $\gtM$. The following
lemma describes the main properties of this action:

\proclaim{Lemma} {\rm a)}  The action of $\Gamma $ on $\gtM$ is proper.

{\rm  b)} For every locally compact space $Z$ endowed with a
proper action of $\Gamma ${\rm ,} there is a continuous equivariant map{\rm ,}
unique up to equivariant homotopy{\rm ,} $Z\ra \gtM$.
\endproclaim

\demo{Proof} a) For every continuous function with compact support
$\varphi$ on $X$, such that $0\le \varphi\le1$, let
$U_\varphi $ denote the set of measures $\lambda\in \gtM$ such
that $\lambda(\varphi )>1/2$. Clearly the sets $U_\varphi $ form
an open covering of $\gtM$. Moreover, if $\varphi $ and
$\psi$ have disjoint supports, $U_\varphi $ and
$U_\psi$ are disjoint; hence, for every continuous function
$\varphi$ with compact support $K\subset X$, the set
$$\{\,g\in \Gamma \,, \;gU_\varphi \cap U_\varphi\not
=\emptyset\,\}\subset \{\,g\in \Gamma \,, \;g(K) \cap
K\not=\emptyset\,\}$$ is relatively compact in $\Gamma $.

  b) Since $\gtM$ is a convex set, any two maps $Z\ra \gtM$ can be
joined by a linear homotopy. This proves uniqueness (up to homotopy).

Let us prove existence. First, assume that $ X=\Gamma $ with the
action by left translations. Let $c$ be a positive continuous cut-off function on $Z$.
This means, by definition, that the support of $c$ has compact intersection
with the saturation of any compact subset of $Z$ and, for every $z\in Z$,
$\int_{\Gamma } c(g^{-1}z)dg=1$. For any $z\in Z$, consider the function
on $\Gamma $: $g\mapsto c(g^{-1}z)$. The product of this function with the
Haar measure on $\Gamma $ is a probability measure on $\Gamma $. The
map $Z\lra \gtM$ associating to $z$ this measure is equivariant.

In general, choosing $x\in X$, we get an equivariant map $\Gamma\to X: g\mt gx$; the
corresponding map on measures is an equivariant map from the space
of measures on $\Gamma $ to the corresponding space of measures on $X$.
\enddemo

It follows from Lemma 4.1 that the space $\gtM$ associated with any
proper $\Gamma $-space $X$ is equivariantly homotopy equivalent to the
universal $\Gamma $-space $\beg $. However, the space $\gtM$ is
too big. We prefer to deal with some subspaces of this space.

For this, assume moreover that $X$ is endowed with a
$\Gamma $-invariant metric. For $k\in{\bf R}_+$, let
$\gtM_k\subset \gtM$ denote the set of probability measures on
$X$ whose support has diameter $\le k$. Note that, if every
bounded set of $X$ is relatively compact, then for every $k\in{\bf
R}_+\,,\; \gtM_k$ is a closed subset of
$\gtM$, hence locally compact.

Indeed, a positive measure $\mu$ has support of diameter $\le k$ if
and only if
$ \mu(f)\mu(g)=0$ for every pair of functions $f,g\in C_c(X)$ such
that the distance between
their supports is $>k$. Therefore, the set ${\frak N}_k\subset \gtM$ of
measures of support of
diameter $\le k$ is a closed subset of $\gtM$.
For any continuous function with compact support
$\varphi$ on $X$, let
$U_\varphi $ denote the set of positive measures $\lambda\in \gtM$
such that $\lambda(\varphi )>1/2$.
Let $0\le \varphi\le1$. If  every
bounded set of $X$ is relatively compact, there exists a $\psi\in
C_c(X)$ such that $0\le \psi\le 1$ and
$\psi(x)=1$, for every $x\in X$ with distance $\le k$ to the
support of $\varphi$.
Then, for $\mu\in U_\varphi\cap {\frak N}_k$ we have
$\|\mu\|=\mu(\psi)$. Since the sets
$U_\varphi\cap {\frak N}_k$ form an open covering of ${\frak N}_k$, the set
$\gtM_k$ of probability measures in ${\frak N}_k$ is a closed
subset\break of ${\frak N}_k$.

For any locally compact space $Z$ endowed with a proper action of
$\Gamma $, such that the quotient $Z/\Gamma $ is compact, there
exists a $k\in{\bf R}_+$ and a continuous equivariant map
$Z\ra \gtM_k$. Moreover, if $f_0$ and
$f_1$ are two such maps, they are homotopic in some $
\gtM_N$ for $N\ge k$.

Let then $\gtM'$ be the telescope of the spaces $\gtM_k$. Let $Z$
be a locally compact, $\sigma$-compact space endowed with a proper
action of
$\Gamma $. Choose a proper function $\varphi:Z/\Gamma \ra {\bf
R}_+$.
There exist:
\begin{itemize}
\item[--] an increasing sequence $k_n\in \R_+$, 
\item[--] an equivariant map $f_n:\varphi^{-1}([0,n])\ra \gtM_{k_n}$,
\item[--]  a sequence $N_n$ with $N_n\ge k_{n+1}$,
\item[--]  an equivariant homotopy $F_n:\varphi^{-1}([0,n])\times [0,1]\ra \gtM_{N_n}$
joining $f_n$ and the restriction of $f_{n+1}$.
\end{itemize}

Let $\psi:\R_+\ra \R_+$ be a continuous increasing function such that\break
$\psi(n)\ge N_n$.
Set then $f(x)=(F_n(x,\varphi(x)-n+1),\psi\circ \phi(x))$ if $\varphi(x)\in
[n-1,n]$.
This is a continuous equivariant map $f:Z\ra \gtM'$.

Moreover, one
may use the same construction for homotopies. It follows that
$\gtM'$ satisfies the conclusion of Lemma 4.1.b) for $\gtM$.
Therefore, the spaces $\gtM$ and
$\gtM'$ are $\Gamma $-equivariantly homotopy equivalent.

For us, it will be sufficient to think of $\gtM$ as of an inductive limit (in
the sense of homotopy theory) of spaces $\gtM_k$.

Assume,   furthermore, that our space $X$ has bounded coarse
geometry. Let $\mu$ be a $\Gamma$-invariant measure on $X$ such that
for any $x\in X$, $\mu (B(x,\delta ))\ge 1$ and  for any $ R>0$
there exists $K(R)>0$ such that for any  subset $S\subset X$ of diameter
$\le R$, $\mu (S)\le K(R)$ (see Proposition 3.2).

\numbereddemo{Definition} For any $N\in \R_+$, define a linear map $\tau :\gtM_N \to
L^2(X;\mu)$ by the formula: $\tau (\nu)=\int_{X}
\chi_{B(x,\delta)}\,d\nu(x)$, where
$\chi_Z$ is the characteristic function of the set $Z$ in $X$. 
\enddemo

\proclaim{Lemma} {\rm a)} Let $R\in \R_+$ and $g$ be a bounded
$\mu$\/{\rm -}\/measurable function on
$X$ such that the diameter of its support is $\le R$. Then
$\|g\|_1 K(R)^{-1/2}\le
\|g\|_2\le \|g\|_\infty K(R)^{1/2}$.
\vglue4pt
{\rm  b)} The image $\tau(\gtM_N)$ in $L^2(X;\mu)$ is contained
between the spheres of radii $K(N+2\delta)^{-1/2}$ and
$K(N+2\delta)^{1/2}$.
\vglue4pt
 {\rm c)} If $X$ is locally finite{\rm ,} the map $\tau
:\gtM_N\ra L^2(X;\mu)-\{0\}$ is continuous and proper in the topology induced
by the weak topology of $L^2(X;\mu)$. Therefore $\tau (\gtM_N)$ is a locally
compact proper
$\Gamma$\/{\rm -}\/space.
\endproclaim

\demo{Proof} a) Let $\chi$ be the characteristic function of the support of
$g$. Replacing $\mu$ by $\chi\cdot \mu$ does not change the
$p$-norms of $g$. Now the total mass of $\chi\cdot \mu$ is
$\le K(R)$, and a) follows.

\vglue12pt b) Let $\nu \in \gtM_N$. As the total mass of $\nu$ is $1$, we
deduce that
$\|\tau(\nu)\|_\infty\le 1$. Since the $\mu$-measure of any open ball of
radius $\delta$ is $\ge1$, we have: $\|\tau(\nu)\|_1\ge 1$. Now b) follows
from a) because the support of $\tau(\nu)$ has diameter $\le N+2\delta$.

\vglue12pt c) The continuity of $\tau$ is obvious since $X$ is discrete. If
$\mu _n$ is a sequence converging to the point at infinity of the one point
compactification of $\gtM_N$, its support goes to infinity in $X$,
and so does the
support of $\tau(\mu _n)$. As $\|\tau (\mu_n)\|$ is bounded by b),
$\tau(\mu _n)$
converges weakly to $0$.
\enddemo
\vglue4pt
 {\it Remark}. In the case of a non locally finite $X$,
assertion  c) remains true if we replace $\chi_{B(x,\delta)}$ by a
continuous approximation.
\vglue16pt

When the space $X$  is locally finite, each $\gtM_k$ is a locally finite
simplicial complex, called a {\it Rips complex}. Therefore $\gtM$ may be
considered as an inductive limit (in the sense of homotopy theory) of
Rips complexes $\gtM_k$.

We remark here that such simplicial presentation of $\gtM$ exists for
any countable discrete group: we may take $X=\Gamma $ and define the
distance by means of a proper length function $\ell$; for example: let
$(g_n)_{n\in{\bf N}}$ be a set of generators for $\Gamma $ and let
$\ell(g)$ be the minimum of $\dst\sum_{i=1}^p\,|r_i|\,(n_i+1)$
over all decompositions $g=g_{n_1}^{r_1}\ldots g_{n_p}^{r_p}$.

\demo{{R}emarks}
a) For $r\in [0,1)$,
the space of finite positive measures on $X$ with total mass contained
in $(r,1]$ is locally compact, but the action of $\Gamma $ on this
space is proper if and only if $r\ge 1/2$.
\vglue12pt  b) Assume that $X$ is endowed with a
$\Gamma $-invariant measure $\mu$. Another realization of the
classifying space for proper actions is the set of nonnegative
$L^2$-functions of norm in the interval $(2^{-1/2},1]$.
\enddemo

\def\RK{RK}
 \vglue-8pt
\section{Novikov's conjecture: an outline of our approach}
 
Let $\Gamma $ be a countable discrete group. There are several
conjectures associated with the Novikov conjecture for $\Gamma $ (see
[K2, 6.4]). All these conjectures deal with the classifying space for free
proper actions of $\Gamma $, usually denoted by $E\Gamma $. The so-called Strong Novikov Conjecture is the
statement that a natural homomorphism
$\beta:RK_*^{\Gamma }(E\Gamma )=RK_*(B\Gamma )\ra
K_*(C^*(\Gamma ))$ is rationally injective. It is known that \pagebreak this statement
implies the Novikov conjecture for $\Gamma $. However, we prefer to deal with
the universal space for proper actions $\beg $ instead of $E\Gamma $. In view
of the discussion of the previous section, we can consider $\beg$ as a
locally compact space.

As  explained in [BCH], the group $RK_*^{\Gamma }(E\Gamma )\otimes
{\bf Q}$ is a subgroup of $RK_*^{\Gamma }(\beg)\break\otimes {\bf Q}$.
Also in [BCH], there is defined a natural homomorphism
$RK_*^{\Gamma }(\beg)\ra K_*(C^*(\Gamma ))$, which we
still prefer to call $\beta$ (we define this map below), and which rationally
coincides on $RK_*^{\Gamma }(E\Gamma )$ with the above
homomorphism $\beta$.

 Let us fix some notation related with crossed products. Let
$\Gamma $ be a locally compact group acting (on the left) on a
$C^*$-algebra $B$. Denote by $dg$ the left Haar measure of
$\Gamma $. The algebra $B$ is contained in the multiplier
algebra of the crossed product $C^*(\Gamma ,B)$ and there is a
canonical strictly continuous morphism $g\mt u_g$ from
$\Gamma $ to the unitary group of the multiplier algebra of the
crossed product $C^*(\Gamma ,B)$. For $b\in B$ and $g\in
\Gamma $, we have $u_gbu_g^*=g\cdot b$; moreover, if $F\in
C_c(\Gamma ,B)$, the multiplier $\int F(g)u_g\,dg$ is actually an
element of $C^*(\Gamma ,B)$, and these elements form a dense
subalgebra of $C^*(\Gamma ,B)$.

Let $\Gamma $ act properly (on the left) on a locally compact
space $Y$. If the action of $\Gamma $ is free, the
algebras $C_0(Y/\Gamma )$ and $C^*(\Gamma ,C_0(Y))$ are Morita
equivalent. In general, we have only a Hilbert $C^*(\Gamma ,C_0(Y))
$-module $E_Y$ and an isomorphism between $C_0(Y/\Gamma )$ and
$\K(E_Y)$ (which is enough for our purposes). To define $E_Y$,
consider $C_c(Y)$ as a left $\Gamma
$-module. For any $h,h_1,h_2\in
C_c(Y)$ and $f\in C_c(\Gamma ,C_0(Y))$, put 
\begin{eqnarray*}
h\cdot
f&=&\int_{\Gamma }g(h)\cdot g(f(g^{-1}))\cdot\nu (g)^{-1/2} dg\in
C_c(Y),\\[6pt]
\langle h_1,h_2\rangle (g)&=&\nu (g)^{-1/2}\overline{h_1}g(h_2)\in
C_c(\Gamma ,C_0(Y))\qquad (g\in\Gamma)\,,
\end{eqnarray*}
 where $\nu (g)$ is the modular
function of
$\Gamma $. One can easily check that $C_c(Y)$ is a submodule of the
pre-Hilbert module $C_c(\Gamma ,C_0(Y))\subset C^*(\Gamma
,C_0(Y))$ (considered as a module over itself). The embedding $i$ is given by
the formula:
$i(h)(g)=\nu (g)^{-1/2}\cdot c^{1/2}\cdot g(h),$ where $ c$ is a positive continuous
cut-off function on $Y$ (this means, by definition, that the support of
$c$ has compact intersection with the saturation of any compact subset of
$Y$ and, for every $y\in Y$, $\int_{\Gamma } c(g^{-1}y)dg=1$). It
follows that the above inner product on $C_c(Y)$ is positive, so we can take
completion which will be denoted by $E_Y$.

One checks immediately that ${\cal K} (E_Y)$ is isomorphic to
$C_0(Y/\Gamma )$ (acting by pointwise multiplication on
$C_c(Y)$).

If $Y/\Gamma $ is compact, then $\K(E_Y)\simeq C(Y/\Gamma)$ is
unital, so $E_Y$ is a finitely generated projective
$C^*(\Gamma,C_0(Y))$-module. Therefore $E_Y$ defines an element of
$K_0(C^*(\Gamma ,C_0(Y))$ which will be denoted by $\lambda _Y$. Let
$f:Y_1\ra Y_2$ be a continuous proper $\Gamma $-map between two proper
locally compact $\Gamma $-spaces with compact quotient. We
obviously have $\lambda_{Y_1}=f^*(\lambda_{Y_2})$ (where
$f^*:K_0(C^*(\Gamma ,C_0(Y_2)))\ra K_0(C^*(\Gamma ,C_0(Y_1)))$ is the
map induced \pagebreak by $f$).

  \centerline{\bf The Baum-Connes map $\beta$}
\vglue6pt

Let $Y$ be a proper locally compact $\Gamma $-space with compact
quotient. Define $\beta_Y:K_\Gamma ^i(C_0(Y))\ra K_i(C^*(\Gamma ))$
by $\beta_Y(x)=\lambda_Y\otimes_{C^*(\Gamma,C_0(Y))} j_\Gamma (x)$. If
$f:Y_1\ra Y_2$ is a continuous proper $\Gamma $-map between two proper
locally compact $\Gamma$-spaces with compact quotient, we
obviously have $\beta_{Y_1}=\beta_{Y_2}\circ f_*$.

\numbereddemo{Definition} Let $\Gamma $ be a locally compact group acting {\it
properly}\/ on a locally compact space $Z$. Put
$RK_i^\Gamma (Z)= \dst\lim_{\lra} K_\Gamma ^i(C_0(Y))$, where the
inductive limit is taken on $Y$ running over $\Gamma $-invariant
closed subsets of $Z$ such that $Y/\Gamma $ is compact. The
Baum-Connes map $\beta :RK_*^\Gamma (\beg)\ra
K_*(C^*(\Gamma ))$ is the map defined at the inductive limit level by
the maps $\beta_Y$. Denote also by $\beta _{\rm red}: RK_*^\Gamma
(\beg)\ra K_*(C_{\rm red}^*(\Gamma ))$ the composition of $\beta$ with
the $K$-theory map associated with the homomorphism $C^*(\Gamma )\ra
C_{\rm red}^*(\Gamma )$.
\enddemo

This map coincides with the map $\mu$ defined in [BCH].

Moreover, if $A$ is a $\Gamma $-algebra, we set $\RK_i^\Gamma
(Z;A)= \dst\lim_{\lra} KK_\Gamma ^i(C_0(Y),A)$, where the inductive
limit is taken on $Y$ running over $\Gamma $-invariant closed
subsets of $Z$ such that $Y/\Gamma $ is compact. One defines in
the same way the Baum-Connes map $\beta ^A:\RK_*^\Gamma (\beg;A)\ra
K_*(C^*(\Gamma ,A))$ and $\beta_{\rm red} ^A:\RK_*^\Gamma (\beg;A)\ra
K_*(C_{\rm red}^*(\Gamma ,A))$.

\vglue2pt
Theorem 1.1 is a consequence of the following theorem, which is the main
result of this paper:

\proclaim{Theorem} For any discrete group $\Gamma $ acting properly by isometries on a
weakly bolic{\rm ,} weakly geodesic metric space of bounded coarse geometry and
every $\Gamma $\/{\rm -}\/algebra $A${\rm ,} the Baum\/{\rm -}\/Connes map $\beta_{\rm red}^A$
is injective.
\endproclaim

It follows that $\beta ^A$ is also injective. We will prove Theorem~5.2
in Sections 7 and 8.

Let $A$ and $B$ be $\Gamma $-algebras. For $x\in
\RK_i^\Gamma (Z;A)$ and $y\in KK^\Gamma(A,B)$, one may form
the $KK$-product $x\otimes _A y\in \RK_i^\Gamma (Z;A)$. One
obviously has:

\proclaim{Proposition}
Let $A$ and $B$ be $\Gamma $\/{\rm -}\/algebras and $a \in
KK^\Gamma (A,B)$. For $x\in \RK_*^\Gamma (\beg;A)${\rm ,}
$\beta ^A(x)\otimes_{C^*(\Gamma ,A)} j_\Gamma (a)=\beta
^B(x\otimes_A a)$. If $\beta ^B$ is an isomorphism and if there
exists an element
$b\in KK^\Gamma (B,A)$ such that $a \otimes_B b=1_A\in
KK^\Gamma (A,A)$ then $\beta ^A$ is an isomorphism. The same
holds if $\beta^A$ and $\beta^B$ are replaced by  $\beta^A_{\rm red}$ and
$\beta^B_{\rm red}$.
\endproclaim

\centerline{\bf Descent isomorphism}
\vglue6pt

Let $\Gamma $ be a locally compact group, $Y$ be a proper $\Gamma $-space, not
necessarily\break $\Gamma$-compact. Denote by $\Lambda_Y$ the element
$$(E_Y,0)
\in {\cal R}KK(Y/\Gamma ; C_0(Y/\Gamma ),C^{*}(\Gamma
,C_0(Y))).$$

\proclaim{Theorem} Let $\Gamma$ be a locally compact group{\rm ,}
$Y$  a proper $\Gamma$\/{\rm -}\/space and $B$ a $\Gamma -C_0(Y)$\/{\rm -}\/algebra. Then{\rm ,}
for $i=0,1${\rm ,} the map $x\mt \Lambda_Y\otimes_{C^*(\Gamma,C_0(Y))}
j_{\Gamma }(x)$ is an isomorphism $${\cal R}KK^i_{\Gamma }(Y;C_0(Y),B)\simeq
{\cal R}KK^i(Y/\Gamma ; C_0(Y/\Gamma ),C^{*}(\Gamma ,B)).$$ If
$Y/\Gamma $ is compact{\rm ,} $${\cal R}KK^i_{\Gamma }(Y;C_0(Y),B)\simeq
K_i(C^{*}(\Gamma ,B)).$$
\endproclaim

Before we give the proof of this theorem, we want to state a result
which will be used in the proof. This is a generalization of the stabilization
theorem for Hilbert modules ([K1]) involving proper group actions. Some
generalizations of this kind are already known (cf.\ [P, 2.9], for example).

\proclaim{Proposition}
Let $\Gamma$ be a locally compact group{\rm ,}
$Y$ a proper $\Gamma$\/{\rm -}\/space and $B$ a $\Gamma -C_0(Y)$\/{\rm -}\/algebra. Assume that
the Hilbert module $\E$ over $B$ is countably generated. Then
$${\cal E}\oplus
(\oplus_1^{\infty}L^2(\Gamma ,B))\simeq\oplus_1^{ \infty}L^2(\Gamma,B).$$
\endproclaim

\demo{Proof} This isomorphism can be obtained in
three steps. First, we embed ${\cal E}$ in
$L^2(\Gamma ,{\cal E})$ as a direct summand using a cut-off function $c$
on $Y$ as follows: $e \mapsto f(g)=g(c)^{1/2}e$. (The projection
$L^2(\Gamma ,{\cal E} )\rightarrow {\cal E}$ is given by
$f\mapsto\int_{\Gamma }f(g)g(c)^{1/2}dg$.) Next, we use the usual infinite
sum trick: $\E\oplus\E^\perp\oplus\E\oplus\E^\perp\oplus...$,
to show that ${\cal E}\oplus (\oplus_1^{\infty}L^2(\Gamma ,{\cal E}
))\simeq\oplus_1^{\infty}L^2(\Gamma  ,{\cal E})$. Finally, we use the
stabilization theorem without group action to get 
\vglue12pt
\hfill ${\displaystyle L^2(\Gamma
,{\cal E})\oplus (\oplus_1^{\infty}L^2(\Gamma ,B))\simeq
\oplus_1^{\infty}L^2(\Gamma ,B).}$ \enddemo
\vglue12pt

{\it Proof of Theorem} 5.4.  Let $({\cal E},T)\in {\cal
R}KK_{\Gamma }(Y;C_0(Y),B)$ and let $C^*(\Gamma ,{\cal E})$ be the Hilbert
module over
$C^{*}(\Gamma ,B)$ defined in [K2, 3.8]
(in fact, $C^*(\Gamma ,{\cal E})={\cal E}\otimes _B
C^{*}(\Gamma ,B)$). Define the Hilbert module
$\widetilde {{\cal E}}$ over $C^{*}(\Gamma ,B)$ by setting $$\widetilde {{\cal
E}}=E_Y\otimes_{C^*(\Gamma ,C_0(Y))} C^*(\Gamma ,{\cal E}).$$ 

The Hilbert module $\widetilde\E$ can also be constructed as follows. Let $
{\cal E}_c=\break C_c(Y)\cdot {\cal E}$. For any
$e,e_1,e_2\in {\cal E}_c$  and $f\in C_c(\Gamma ,B)$, put
\begin{eqnarray*}
e\cdot
f&=&\int_{\Gamma}g(e)\cdot  g(f(g^{-1}))\cdot\nu (g)^{-1/2} dg\in {\cal
E}_c,\\
\langle e_1,e_2\rangle (g)&=&\nu  (g)^{-1/2}(e_1,g(e_2))_{\cal E}\in
C_c(\Gamma ,B),\pagebreak
\end{eqnarray*}
 where $\nu (g)$ is  the modular function of
$\Gamma$. There is  a natural map of the algebraic tensor product
$C_c(Y)\otimes C_c(\Gamma,\E)$ to $\E_c$ given by
$$f\otimes e\mapsto \int_\Gamma \nu(s)^{-1/2}s^{-1}(f)s^{-1}(e(s))ds$$
which preserves the inner products and the right actions of $C_c(\Gamma,B)$.
This map extends to an isomorphism of $\widetilde\E$ with the completion of
$\E_c$.

An easy argument shows that ${\cal L}(\widetilde {{\cal E}})$ is isomorphic
to the $\Gamma $-invariant part of ${\cal L}({\cal E})$ and that ${\cal K}
(\widetilde {{\cal E}})$ is isomorphic to ${\cal K}({\cal E})^{\Gamma }$ (see
[K2, Def.~3.2]). This means that we can consider $\widetilde
{T}=\int_{\Gamma }g(cT)dg$ as an operator on $\widetilde {{\cal E}}$ (where
$c$ is a cut-off function). The map $ ({\cal E},T)\mapsto (\widetilde
{{\cal
E}},\widetilde {T})$ gives a  homomorphism of ${\cal R}KK^i_{\Gamma
}(Y;C_0(Y),B)$ to
${\cal R}KK^i(Y/\Gamma ; C_0(Y/\Gamma ),   C^{*}(\Gamma ,B))$ which
coincides with the homomorphism $x\mt \Lambda_Y\otimes_{C^*(\Gamma,C_0(Y))}
j_{\Gamma }(x)$.

To prove that this is an isomorphism we apply Proposition 5.5 which allows us
to assume that our initial Hilbert $B$-module ${\cal E}$ is isomorphic to
$\oplus_1^{
\infty}L^2(\Gamma  ,B)$. To finish the proof, it is enough to show that in this
case, $\widetilde {{\cal E}}\simeq\oplus_1^{\infty}C^{*}(\Gamma ,B)$ as a
Hilbert
module over $ C^{*}(\Gamma ,B)$. Of course, we will take only one copy
of $L^2(\Gamma ,B)$ and prove that if ${\cal E}\simeq L^2(\Gamma ,B)$
then $\widetilde {{\cal E}}\simeq C^{ *}(\Gamma ,B)$ as a Hilbert module over
$C^{*}(\Gamma ,B)$. To get this, it will be convenient to consider
$L^2(\Gamma ,B)$ with the right $\Gamma $-action: $g(f)(g_1)=\nu
(g)^{1/2}g(f(g_1g))$, instead of the usual left one. (The two
$\Gamma $-actions, clearly, correspond to each other under the
automorphism $f(g)\mapsto\nu (g)^{-1/2}f(g^{-1})$ of $L^ 2(\Gamma
,B)$.)  With this convention, the desired isomorphism $\widetilde {{\cal
E}}\simeq C^{*}(\Gamma ,B)$ is given by the formula: $\tilde
{e}(g)\mapsto g(\tilde {e}(g))$.
\hfill\qed

\vglue12pt \centerline{\bf Proper algebras}
 
\numbereddemo{Definition} A $\Gamma $-algebra is said to be {\it proper}\/ if
it is a $\Gamma -C_0(Z)$-algebra for some proper $\Gamma $-space $Z$.
\enddemo

Since every proper $\Gamma$-space maps equivariantly to $\beg$, a
$\Gamma $-algebra is  proper if and only if it is a $\Gamma
-C_0(\beg)$-algebra.

The following proposition is a particular case of some results of [Tu,
\S 5]. As
some of the statements and proofs there are a little   too imprecise, we
prefer to give
a complete proof here.

\proclaim{Proposition} Let $\Gamma $ be a second countable locally compact group{\rm ,} $X$ a
second countable locally compact
$\Gamma $\/{\rm -}\/space
and $A$ a nuclear $\Gamma $\/{\rm -}\/algebra. Assume that the $\Gamma${\rm -}algebra
$A\otimes
C_0(X)$ is proper. Then the functor $B\lra RKK^\Gamma (X;A,B)$ is {\rm `}\/half
exact{\rm '. (}\/All algebras are assumed to be separable.\/{\rm )}
\endproclaim

 This means that for every $\Gamma $-equivariant short exact sequence of
$\Gamma$-algebras $$0\to J {\buildrel  i\over\lra }B {\buildrel q\over\lra}
B/J\to 0 ,$$
the sequence $$RKK^\Gamma (X;A,J){\buildrel i_*\over\lra} RKK^\Gamma (X;A,B)
{\buildrel q_*\over\lra} RKK^\Gamma (X;A,B/J)$$ is exact in its middle
term, from which
it follows that we have a six term exact \pagebreak sequence.

 {\it Proof}. We follow the proof of [S, Prop.\ 3.1]. Let us state the intermediate
Lemmas (3.2--3.3 of [S]) in our context.

\proclaim{Lemma} Let $(\E,F)$ be an element in $RKK^\Gamma (X;A,B)$. Put $\overline
\E=\E \,\widehat{\otimes}\,_BB/J$.
 \vglue4pt
{\rm a)} If $q_*(\E,F)$ is degenerate{\rm ,} then $(\E,F)$ is in the image
of $i_*$. 
\vglue4pt
{\rm b)} An operator homotopy $(\overline \E,G_t)$ in $RKK^\Gamma (X;A,B/J)$
with $G_0=F \,\widehat{\otimes}\,1$ can be lifted to an operator homotopy
$(\E,F_t)$ in
$RKK^\Gamma (X;A,B)$ with
$F_0=F$.\endproclaim

The proof of these facts is the same as in the nonequivariant setting:

\demo{Proof}  We have an exact sequence $0\to \K(\E_J)\to \K(\E)\to\K(\overline
\E)\to 0$, where
$\E_J=\{\xi\in \E,\ \langle \xi,\xi\rangle \in J\}$.

\vglue4pt  a) If $q_*(\E,F)$ is degenerate, $(\E_J,F)$ is an element in
$RKK^\Gamma
(X;A,J)$ which, as an element of $RKK^\Gamma (X;A,B)$, is homotopic to $(\E,F)$.

\vglue4pt b) Let ${\cal A}$ (resp.\  ${\cal B}$) be the set of $T\in
\L(\E)$ (resp.\  $T\in \L(\overline \E)$) such that for all
$a\in C_0(X)\otimes A$, the commutator $[a,T]$ is compact and the
function $g\mt
a(g T-T)$ ($g\in \Gamma $) is norm-continuous with compact values. Let also
${\cal I}$ (resp.\  ${\cal J}$) be the set of $T\in \A$
(resp.\ 
$T\in \B$) such that for all $a\in C_0(X)\otimes A$, $Ta$ is
compact.

We claim that the morphism ${\cal A}/{\cal I}\to {\cal B}/{\cal J}$ is
onto. Indeed, let $S\in{\cal B}$. Since the morphism $\hat q:\L(\E)\to
\L(\overline \E)$ is surjective, we can find
$T\in \L(\E)$ with image $S$. Averaging  
$S$ and $T$ with respect to a continuous cut-off function on~$\Gamma$, we may assume that $S$ and $T$ are  
$\Gamma$-continuous (this changes
$S$ by some element of ${\cal J}$). Let $D$ be the (separable) subalgebra of
$\L(\E)$ generated by $\K(\E)$, $C_0(X,A)$ and the translates of $T$ by
$\Gamma$. Set $D_1=D\cap \ker \hat q$. Now thanks to Theorem 1.4 of
[K2], one may
construct a $\Gamma$-continuous, equivariant up to $\K(\E_J)$, element $M\in
\L(\E)$ which
commutes with $A$ and $T$ up to $\K(\E_J)$, such that $0\le M\le 1$,
$M D_1\subset
\K(\E_J)$ and $(1-M) \K(\E)\subset \K(\E_J)$
\footnote{According to [K2], $M$ can be chosen
as an element $M_0$ of $\L(\E_J)$. If ${\cal K}$ is an ideal in a $C^*$-algebra
${\cal D}$, the algebra ${\cal M}({\cal D},{\cal K})$ of multipliers $T$ of
${\cal D}$ such that $T{\cal D}+{\cal D}T\subset \K$ embeds both in ${\cal
M}(\K)$ and ${\cal
M}({\cal D})$; take $M\in \L(\E)$ such that $(1-M)\in {\cal
M}(\K(\E),\K(\E_J))$
with image $1-M_0$ in ${\cal M}(\K(\E_J))=\L(\E_J)$.}.
From the last inclusion, it follows that
$1-M\in \ker \hat q$, whence $S=\hat q(MT)$. Now, the elements $[T,a],
a(gT-T)$ belong to
$D\cap \hat q^{-1}(\K(\overline \E))$. Note that an element
$x\in D\cap \hat q^{-1}(\K(\overline \E))$  can be written as a sum $x=y+z$
where $y\in
\K(\E)$ and $z\in D_1$. Therefore $Mx\in \K(\E)$. It follows easily that
$MT\in \A$.

Let $U$ (resp.\  $V$) denote the set of self-adjoint elements of
degree $1$ and square $1$
in ${\cal A}/{\cal I}$ (resp.\  ${\cal B}/{\cal J}$). The map $U\to
V$ obviously
satisfies the homotopy lifting property. The result follows.
\enddemo

By Lemma 5.8, if $q_*(\E,F)$ is operator homotopic to a degenerate element,
its class is in
the image of $i_*$. Now, if the class of $q_*(\E,F)$ is $0$, there exists a
degenerate
element $(\E',F')$ in $RKK^\Gamma (X;A,B/J)$ such that $q_*(\E,F)\oplus
(\E',F')$
is operator homotopic to a degenerate element. Furthermore, if a degenerate
element $(\E'',F'')$ of $RKK^\Gamma (X;A,B/J)$ contains $(\E',F')$ as a
direct
summand, then obviously $q_*(\E,F)\oplus (\E'',F'')$
is operator homotopic to a degenerate element.

Therefore, to end the proof of our proposition we just need to prove the
following
analogue of Lemma 3.5 in [S]:
\vglue-4pt
\proclaim{Lemma} For every degenerate element $(\E',F')$ in $RKK^\Gamma
(X;A,B/J)${\rm ,} there
exists a degenerate element $(\widetilde \E,\widetilde F)$ in $RKK^\Gamma
(X;A,B)$ such that\break $(\widetilde\E \,\widehat{\otimes}\,_BB/J,\widetilde
F \,\widehat{\otimes}\,1)$ contains $(\E',F')$ as a direct summand.
\endproclaim

\vglue-4pt
{\it Proof}. A representation of $A\otimes C_0(X)$ is just a pair of commuting
representations. Now, since the left and right actions of $C_0(X)$
have to be the
same, the only difference between elements of $KK_\Gamma(A,B\otimes
C_0(X))$ and
$RKK_\Gamma(X;A,B)$ is the compactness requirements. The degenerate elements
are the same. The representation of $A$ together with the element $F'$ define a
representation $A \widehat{\otimes}\,{\cal C}_1\to \L(\E')$. In other words,
degenerate elements in $RKK_\Gamma(X;A,B)$ are just equivariant
$(A\,\widehat{\otimes}\, {\cal C}_1,C_0(X)\otimes B)$-bimodules.

Using an equivariant representation of $A\,\widehat{\otimes}\, {\cal C}_1$ on a
separable Hilbert space $\H$, we may find an equivariant
$(A\,\widehat{\otimes}\,{\cal C}_1,C_0(X)\,\otimes\, B/J)$-bimodule $\E''$
isomorphic to $\H\,\,\widehat{\otimes}\,\,C_0(X)\,\otimes\, B/J  $. Then $\E'$ is a direct
summand in $\E'\oplus \E''$.  By the (nonequivariant) stabilization theorem of
[K1], the $C_0(X)\otimes B/J$-module
$\E'\oplus \E''$ is isomorphic to $\E''$, whence $\L(\E'\oplus \E'')$ is a
quotient of $\L(\H\,\,\widehat{\otimes}\,\,C_0(X)\otimes B)$. Denote by $\pi
:A\,\widehat{\otimes}\,{\cal C}_1\to \L(\E'\oplus \E'')$ the left action. Since
$A\widehat{\otimes}\,{\cal C}_1$ is nuclear, the map $\pi $ admits a completely
positive lifing. Using the Stinespring construction of [K1], we find a Hilbert
$C_0(X)\otimes B$-module $\overline \E$ and a representation\break
$\pi':A \,\widehat{\otimes}\,{\cal C}_1\to \L(\overline \E)$ such that $\pi'\otimes
1=\pi$. Note moreover that the $C_0(X)\break\otimes B$-module $\overline \E$ contains
$\H\,\,\widehat{\otimes}\,\,C_0(X)\otimes B$ as a direct summand, and is therefore
isomorphic to $\H\,\,\widehat{\otimes}\,\,C_0(X)\otimes B$. Consequently,
there exists an
action of $\Gamma
$ on~$\overline \E$.

Note that the action of $A\widehat{\otimes}\,{\cal C}_1$ on $\overline \E$
and the isomorphism $U$ of $\overline \E \,\widehat{\otimes}\,_BB/J$ with
$\E'\oplus \E''$ are not assumed to be $\Gamma $-equivariant.
This is taken care of by tensoring with $L^2(\Gamma )$. Set $\widetilde
\E=L^2(\Gamma)\otimes \overline \E$ as a $C_0(X)\otimes B-\Gamma$-module. The
action $\tilde \pi $ of $A\widehat{\otimes}\,{\cal C}_1$ on
$\widetilde \E$ is given by $(\tilde \pi(a)\xi)(g)=g\cdot
(\pi'(g^{-1}\cdot
a)(g^{-1}\cdot \xi(g))$ ($a\in
A \,\widehat{\otimes}\,{\cal C}_1\,, \ \xi\in \widetilde \E=L^2(\Gamma
,\overline\E)\,,
\;g\in
\Gamma $). It is equivariant. 

We claim that the $(A\,\widehat{\otimes}\,{\cal C}_1,C_0(X)\otimes
B/J)$-bimodules
$\widetilde \E$ and $(\E'\oplus \E'')\otimes L^2(\Gamma) $ are isomorphic.  The
element $\widetilde U\in \L(\widetilde
\E \,\widehat{\otimes}\,_BB/J,
\E'\otimes
L^2(\Gamma ))$ given by
$(\widetilde U\xi)(g)=g\cdot (U(g^{-1}\cdot \xi(g))$ is $\Gamma
$-invariant. Moreover, since the action of $A\,\widehat{\otimes}\,{\cal C}_1$ on
$\E'\oplus \E''$ is
$\Gamma $-equivariant,
$\widetilde U$ intertwines the actions of $A\,\widehat{\otimes}\,{\cal
C}_1$.

We finally prove that the $(A\,\widehat{\otimes}\,{\cal C}_1,C_0(X)\otimes
B/J)$-bimodule $\E'$ is a direct summand of $ \E'\otimes L^2(\Gamma )$.

Let $Y$ be a proper $\Gamma$-space such that $C_0(Y)$ acts in a  nondegenerate
way by central multipliers on $C_0(X)\otimes A$. Let $c:Y\to \C$ be a positive
cut-off function. Let $\Gamma$ act by left translations on $\Gamma$
and diagonally
on $C_0(Y)\otimes L^2(\Gamma)$. Associated to $c$ is an isometry $V_0:C_0(Y)\to
C_0(Y)\otimes L^2(\Gamma)$ given by $V_0(\xi)(y,g)=\xi(y)c (g^{-1}y)^{1/2}$,
where $\xi\in C_0(Y)$ and $V_0(\xi)\in C_0(Y)\otimes L^2(\Gamma)$ is seen as a
function of two variables
$y\in Y$ and $g\in \Gamma$. One checks immediately that $V_0$ is a
$\Gamma$-invariant element of $\L(C_0(Y),C_0(Y)\otimes L^2(\Gamma))$ and
$V_0^*V_0=1$. Now, write
$$C_0(X)\otimes A\,\widehat{\otimes}\,{\cal C}_1=C_0(Y)\otimes
_{C_0(Y)}(C_0(X)\otimes A\,\widehat{\otimes}\,{\cal C}_1)$$ and 
$$C_0(X)\otimes
A\,\widehat{\otimes}\,{\cal C}_1 \otimes L^2(\Gamma)=(C_0(Y)\otimes
L^2(\Gamma)) \otimes
_{C_0(Y)}(C_0(X)\otimes A\,\widehat{\otimes}\,{\cal C}_1);$$  let $$V\in
\L(C_0(X)\otimes A\,\widehat{\otimes}\,{\cal C}_1,C_0(X)\otimes
A\,\widehat{\otimes}\,{\cal C}_1\otimes L^2(\Gamma))$$ be
$V_0\otimes 1$. Since the action of $C_0(Y) $ is central, $V$ intertwines the
natural left actions of $A\,\widehat{\otimes}\,{\cal C}_1$.

It follows that the equivariant
$(A\,\widehat{\otimes}\,{\cal C}_1,C_0(X)\otimes B)$-bimodule
$\E'$ is a direct summand of
$(A\,\widehat{\otimes}\, {\cal C}_1\otimes L^2(\Gamma )) \,\widehat{\otimes}\,_{A\,\widehat{\otimes}\, {\cal C}_1}\E'\simeq \E'\otimes L^2(\Gamma )$
and therefore a direct summand of $(\E'\oplus \E'')\otimes
L^2(\Gamma)\simeq \widetilde \E\otimes _{C_0(X)\otimes
B}(C_0(X)\otimes B/J)$. This
ends the proof.
\hfill\qed

\numbereddemo{{R}emark} Let $\Gamma $ be a locally compact group, $X$ a locally compact
$\Gamma $-space
and $A,A'$ nuclear $\Gamma $-algebras. Assume that the
$\Gamma$-algebras $A\otimes
C_0(X)$ and $A'\otimes C_0(X)$ are proper. Let $0\to J \to B{\buildrel q\over
\lra} B/J\to0$ be a short
exact sequence of
$\Gamma$-algebras and $u$ be an element in $RKK_\Gamma(X;A,A')$. Denote by
$\partial:RKK_\Gamma(X;A,B/J)\to RKK_\Gamma^1 (X;A,J)$ and $\partial
':RKK_\Gamma(X;A',B/J)\to
RKK_\Gamma^1 (X;A',J)$ the connecting maps associated with the exact
sequences. These connecting maps
are obtained by composing the map $B(0,1)\to C_q$ and the inverse of the
map $e:J\to C_q$ where
$C_q=B[0,1)/J(0,1)$ is the cone of $q$. Therefore, for any $x\in
RKK_\Gamma(X;A',B/J)$
we have $\partial (u\otimes_{A'} x)=u\otimes_{A'}\partial '( x)$.
\enddemo

Using now Corollary A.4 of the appendix, for any $\Gamma $-invariant
closed subset $Y$ of $\beg$ and any $\Gamma $-algebra $B$,
we obtain an isomorphism: $KK_\Gamma ^i(C_0(Y),B)=E_\Gamma ^i(C_0(Y),B)$,
and therefore $\RK_i^\Gamma (\beg;B)$ is equal to the group $ E_\Gamma
^i(\beg,B)$ (of [GHT]). Moreover, for any proper algebra $B$,
$C^*(\Gamma ,B)=C_{\rm red}^*(\Gamma ,B)$. This allows us to apply certain
methods and results of [GHT] to $KK$-theory. In particular we obtain

\proclaimtitle{cf.\ [GHT, Th. 13.1]}
\proclaim{Proposition}
 Assume that the $\Gamma $\/{\rm -}\/algebra $B$ is
proper. Then the Baum\/{\rm -}\/Connes
homomorphisms $\beta ^B$ and $\beta_{\rm red} ^B$ are split surjective.
If the group $\Gamma $ is
discrete{\rm ,} these homomorphisms are isomorphisms.
\endproclaim

\demo{Sketch of proof} Let us describe the inverse map:

Note that $K_i(C^*(\Gamma ,B))$ is the inductive limit of
$K_i(C^*(\Gamma ,C_0(U)B))$ where $U$ runs over open $\Gamma
$-invariant subsets of $\beg $ such that $U/\Gamma $ is
relatively compact.

Let $U\subset Y\subset \beg$ be $\Gamma $-invariant
subsets of $\beg$ with $U$ open and $Y/\Gamma $
compact. Theorem 5.4 gives an isomorphism
$$K_i(C^*(\Gamma ,C_0(U)B))\simeq {\cal R}KK^i_{\Gamma }(Y;C_0(Y),C_0(U)B).$$
Denote by $\alpha _{U,Y}$ the composition
$$K_i(C^*(\Gamma ,C_0(U)B))\simeq {\cal
R}KK^i_{\Gamma }(Y;C_0(Y),C_0(U)B)\rightarrow
KK^i_{\Gamma }(C_0(Y),C_0(U)B).$$

Let $U\subset V\subset Y\subset Z$ be $\Gamma $-invariant
subsets of $\beg$ with $U,V$ open and $Y/\Gamma ,
Z/\Gamma $ compact. Denote by $u:C_0(U)B\rightarrow C_0(V)B$ the
natural inclusion map and $q:C_0(Z)\rightarrow C_0(Y)$ the restriction
map. We obviously have $q^*\circ \alpha _{U,Y}=\alpha _{U,Z}$ and
$u_*\circ \alpha _{U,Y}=\alpha _{V,Y}\circ u_*$, from which it follows
that we may take an inductive limit in
$Y$ and get a map
$\alpha _U :K_i(C^*(\Gamma ,C_0(U)B))\rightarrow \RK_i^\Gamma
(\beg ,B)$. Moreover, since $\alpha _U=\alpha _V\circ u_*$,
the maps $\alpha_U$ define a morphism
$\alpha :K_i(C^*(\Gamma ,B))\ra \RK_i^\Gamma (\beg,B)$.
\vglue2pt
It is easy to see that $\beta ^B\circ
\alpha$ is the identity of $K_i(C^*(\Gamma ,B))$.
But the fact that the composition $\alpha \circ
\beta ^B$ is also the identity in the case of a discrete group $\Gamma$ is
more complicated (cf.\ [GHT]).
\enddemo
\vglue4pt

\centerline{\bf Sufficient conditions for the injectivity of the Baum-Connes
map}
\vglue12pt

To establish the injectivity of the Baum-Connes map, we will use the
following simple result, in which $\widehat\otimes$ stands for (graded)
minimal or maximal tensor products:

\proclaim{Lemma} Let $\Gamma $ be a locally compact group and $B$ a
$\Gamma $\/{\rm -}\/algebra. Assume that for every closed $\Gamma $\/{\rm -}\/invariant
subset $Y\subset \beg$ with compact quotient there exist a
$\Gamma $\/{\rm -}\/algebra $A$ and elements $\eta\in KK^\Gamma_i (\C,A)$
and $d\in KK^\Gamma_i(A,\C)${\rm ,} such that $\beta ^{A \,\widehat{\otimes}\,B}$ {\rm (}\/resp.\
$\beta_{\rm red} ^{A \,\widehat{\otimes}\,B}${\rm )} is injective and
$p_Y^*(\eta\otimes_A d)=1_Y${\rm ,}
where $p_Y$ is the map $Y\ra {\rm point}$ and $p_Y^*$ is the map
$KK^\Gamma(\C,\C)\to RK^0_\Gamma (Y)$.
Then the Baum\/{\rm -}\/Connes map
$\beta ^B$ {\rm (}resp.\ $\beta _{\rm red}^B${\rm )} is injective.
\endproclaim

\demo{Proof} Indeed, let $z\in \ker\beta ^B$; there exists $Y$ and there
is a $$y\in KK_*^\Gamma (C_0(Y),B)$$ with image $z$. Take $A,\eta ,d$
corresponding to $Y$. As $z$ is in the image of $KK_*^\Gamma
(C_0(Y),B)$, we have $z=z\otimes_\C \eta\otimes_A d$. Then $\beta
^{A \,\widehat{\otimes}\,B}(z\otimes _\C\eta)=\beta ^B (z)\otimes _{C^*(\Gamma ,B)}
j_\Gamma (\sigma _B(\eta))=0$, hence $z\otimes_\C \eta=0$ and
$z=z\otimes_\C \eta\otimes_A d=0$.
\enddemo

Combining Lemma~5.12 and Proposition~5.11, we get:

\proclaim{Proposition}
Assume that the group $\Gamma $ is discrete{\rm ,} and for every
closed $\Gamma $\/{\rm -}\/invariant subset $Y\subset \beg$ with compact
quotient there exist
\begin{itemize}
\ritem{--} a proper $\Gamma $-algebra $\A${\rm ;}  
\ritem{--} elements $\eta\in KK_i^\Gamma (\C,\A)$ and $d\in KK_i^\Gamma
(\A,\C)$ such that $p_Y^*(\eta\otimes_\A d)\break =1_Y${\rm ,} where $p_Y$ is the
map $Y\ra {\rm
point}$, $p_Y^*$ is the map $KK^\Gamma(\C,\C)\to RK^0_\Gamma (Y)$.
\end{itemize}
Then $\beta_{\rm red} ^B$ is injective for every
$\Gamma $-algebra $B$.
\endproclaim

\demo{Proof} Let $B$ be a $\Gamma $-algebra. The
algebra $\A \,\widehat{\otimes}\,B$ is proper. By Proposition~5.11, $\beta
_{\rm red}^{\A \,\widehat{\otimes}\,B}$ is an isomorphism. Now the assertion follows from Lemma~5.12.
\phantom{noice}
\enddemo

Assume now that our discrete group $\Gamma $ acts properly by isometries on a
locally finite space $X$ of bounded coarse geometry. As it was explained
in the previous section, any proper $\Gamma $-space $Y$, such that
$Y/\Gamma $ is compact, admits a $\Gamma $-equivariant map into some Rips
complex $\gtM_k$. As an immediate corollary of the previous proposition we
get:

\proclaim{{C}orollary} Assume that for every $k\in {\bf R}_+${\rm ,} there exist
\begin{itemize}
\ritem{--} a proper $\Gamma $-algebra $\A_k${\rm ;}
\ritem{--} elements $\eta_k\in KK^\Gamma_i
(\C,\A_k)$ and
$d_k\in KK^\Gamma_i (\A_k,\C)$ such that\break $p_{\gtM_k}^*(\eta
_k\otimes_{\A_k} d_k)=1_{\gtM_k}${\rm ,} where
$p_{\gtM_k}$ is the map ${\gtM_k}\ra {\rm point}${\rm ,}
$p^*_{\gtM_k}$ is the map $KK^\Gamma(\C,\C)\to RK^0_\Gamma
(\gtM_k)$.
\end{itemize}
  Then $\beta_{\rm red} ^B$ is injective for every
$\Gamma $-algebra $B$.\hfill\qed
\endproclaim
\vglue-14pt
 
\section{The $\gamma$ element}
\vglue-4pt
This section contains one of the main ingredients of the proof of
Theorem~5.2. Namely, assuming that $(X,d)$ is a proper
$\Gamma $-space which is locally finite, weakly geodesic, has bounded coarse
geometry and satisfies a condition somewhat weaker than weak bolicity, we
construct an element
$\gamma_k\in KK_{\Gamma }({\bf C},{\bf C})$ such that
$q^*(\gamma_k)=1_{\gtM_k}\in RK^0_{\Gamma }(\gtM_k)$ (where
$q$ is the projection $\gtM_k\ra\pt$).

We fix a metric space $(X,d)$. Here is some additional notation that we
will use:
 
For $N\in{\bf R}_+$, let $\Delta_N$ denote the set of all nonempty
finite subsets
of $X$ of diameter $\le N$. (Clearly, if $X$ is locally finite,
$\Delta_N$ is a combinatorial complex, the geometric realization of
which is $\gtM_N$.)

For any
$S\in \Delta_N$, set $U_S=\dst\bigcap_{y\in S} \overline B(y,N)=\{\,z\in
X\,,\;S\cup \{z\}\in \Delta_N\,\}$.

We begin with the following:

\proclaim{Lemma} Assume $(X,d)$ is weakly $\delta$\/{\rm -}\/geodesic. Let
$x\in X$ and $S\in \Delta_N$ be such that $ x\not\in U_S$.
For all $z\in U_S${\rm ,} $\sup\{\, d(z,y)\,,\; y\in S\,\}\ge
N+d(x,U_S)-d(x,z)-2\delta$.
\endproclaim

\demo{Proof} Let $\eta \in \R$ be such that $d(x,z) - d(x,U_S) <\eta$;
we must prove that there exists $c\in S$ such that $d(c,y)>
N-2\delta-\eta$. As $(X,d)$ is weakly $\delta$-geodesic, there
exists a point
$b\in X$ such that $d(z,b)\le
\eta+2\delta$ and $d(x,b)\le d(x,z)-\eta$. Since $d(x,b)\le
d(x,z)-\eta<d(x,U_S)$, it follows that $b\not\in U_S$; therefore
there exists $ c\in S\,,\;d(b,c)>N$ whence $d(z,c)\ge d(b,c)-d(z,b)>
N-2\delta-\eta$.
\enddemo

We now fix nonnegative real numbers $\delta,k,N$ such that
$N\ge 8k+22\delta $ and set $\Delta=\Delta_N$. In the sequel of
this section we
assume that $(X,d)$ is weakly $\delta$-geodesic and satisfies the
following
condition (which is a consequence of condition (B2$'$) of weak bolicity):

\begin{itemize}
\item[(C2)] There exists a map $m:X\times X\ra X$ such
that if $x,y,z$ are points of $X\,,$ then $m(x,y)$ is a
$\delta$-middle point of $x,y$ and $d(m(x,y),z)\le \max
(d(x,z),\break d(y,z)) + 2\delta$. If moreover, $d(x,z)\le N$, $d(y,z)\le
N$ and $d(x,y)> N$ then $d(m(x,y),z)< N-4k-10\delta $.
\end{itemize}

\proclaim{Lemma} Let $x\in X$ and $S\in \Delta $. The diameter of $\{\, z\in
U_S\,,\;d(x,z) \le d(x,U_S)+(4k+6\delta )\,\}$ is $\le N$.
\endproclaim

\demo{Proof} If $x\in U_S$ the assertion is obvious since $N\ge
8k+12\delta$. Let $y,z\in U_S$ be such that $d(x,y) \le
d(x,U_S)+4k+6\delta$ and $d(x,z) \le d(x,U_S)+4k+6\delta $;
assume $d(y,z)> N$. By condition (C2) $d(x,m(y,z)) \le d(x,U_S)+
4k+8\delta
$, and for every $c\in S$ we have $d(c,m(y,z))<N-4k-10\delta $,
which is in contradiction with Lemma~6.1.
\enddemo

  From now on, assume that $(X,d)$ is locally finite.

\proclaim{Lemma} Let $x\in X$ and $S,T\in \Delta $. Assume that for every
$a$ in the symmetric difference of $S$ and $T$ we have
$d(a,x)\le d(x,U_S)+4k+6\delta$. Then
\begin{itemize}
\ritem{a)} $d(x,U_T)=d(x,U_S)$. 
 \ritem{b)} For any $b$ in the symmetric difference of $U_S$ and
$U_T${\rm ,} $d(x,b)\ge d(x,U_S)+4k+6\delta$.
\end{itemize}

\endproclaim

\demo{Proof} Using induction on the cardinality of the symmetric difference of
$S$ and $T$ we may assume that this symmetric difference
consists of exactly one element $a$.

Assume first that $a\in T$. Then $U_T\subset U_S$; moreover,
since $T=S\cup \{a\}\in\Delta$, $a\in U_S$. By Lemma~6.2, the
set
$\{\,z\in U_S\,,\; d(x,z)\le d(x,U_S)+4k+6\delta\,\}$ has diameter
less than $N$ and by our assumption it
contains $a$. It is therefore contained in
$U_T$. We have proved that $$\{\,z\in U_S\,,\; d(x,z)\le
d(x,U_S)+4k+6\delta\,\}
\subset U_T\subset U_S;$$ a) and b) follow immediately.

Assume next that $a\in S$. Then $U_S\subset U_T$.

Suppose that $d(x,U_T)<d(x,U_S)$; set
$F=\{\,b\in U_T\,,\; d(x,b)<d(x,U_S)\,\}$ and let $b\in F$ be such
that $d(a,b) =d(a,F)$. As $d(x,b)<d(x,U_S)\,,\; b\not\in U_S$, so
$d(a,b)>N$. Set $b_1=m(a,b)$; by condition (C2), $d(b_1,x) \le
d(a,x)+2\delta$, and there exists a positive real number $
\varepsilon $ such that for all $y\in T\,, \;d(y,b_1) + \varepsilon
<N-4k-10\delta $; we may moreover assume that $2\varepsilon +N<
d(a,b)$. Let $c\in X$ be such that $d(b_1,c)\le 4k+10\delta
+\varepsilon $ and $$d(x,c)\le d(x,b_1)-4k-8\delta -\varepsilon\le
d(x,a)-4k-6\delta-\varepsilon\le d(x,U_S)-\varepsilon
.$$ Then, $c\in F$ and  therefore $$d(a,b)\le d(a,c)\le d(a,b_1)
+4k+10\delta+\varepsilon \le \dst{d(a,b)\over 2}+4k+11\delta
+\varepsilon$$ and $N< d(a,b)-2\varepsilon \le 8k+22\delta$ which
contradicts our hypothesis.

Now a) is proved; we may therefore exchange the roles of $S$ and
$T$;\break b) follows.
\enddemo

\proclaim{Lemma} Let $x\in X$ and $S\in \Delta$ satisfy $\sup\{\,
d(x,y)\,,\; y\in S\,\} > 4k+6\delta $. Then $\sup\{\, d(x,y)\,,\; y\in
S\,\} > d(x,U_S) + 4k+6\delta $. 
\endproclaim

\demo{Proof} If $S$ has one point, the assertion is true since $(X,d)$ is
weakly geodesic. Assume that $\sup\{\, d(x,y)\,,\; y\in S\,\}\le
d(x,U_S)+4k+6\delta $. Let $T$ be a set consisting of one point in
$S$; we get a contradiction using Lemma~6.3.
\enddemo

{\it Notation}. For $R\in {\bf R}_+$, let $I(R)$ be the set
of real numbers
$r\in \R_+$ such that for every quadruple
$x,y,a,b$ of points of $X$ satisfying $d(x,a)+d(y,b)\ge
2R-r\,,\;d(x,y)\le r$ and $ d(a,b)\le 2N$, one has: $d(y,a)+d(x,b)\le
d(x,a)+d(y,b)+2k$. Note that $I(R)$ is an interval in $\R_+$
containing $0$.

\proclaim{Lemma} {\rm a)} For all $R\in\R_+$, $k\in I(R)$.
\vglue4pt
{\rm b)} If $R\le R'\in {\bf R}_+${\rm ,} $I(R)\subset
I(R')${\rm ;} if
$r\in \R_+$ satisfies $r+2(R'-R)\in I(R')${\rm ,} then $r\in
I(R)$.

\vglue4pt
{\rm c)} If the diameter of $X$ is infinite then
$\sup (I(R))\le \sup\{R-N,0\}+k+6\delta\le \sup\{R,6\delta\}+k$.
\endproclaim

\demo{Proof} a) From the inequalities $d(x,b) \le d(x,y)+d(y,b)$ and $d(y,a) \le
d(x,y)+d(x,a)$ we get $d(y,a)+d(x,b)\le d(x,a)+d(y,b)+2d(x,y) $, from
which the first assertion follows.

\vglue4pt In b), the first assertion is obvious. To prove the second
assertion, set\break $r'=r+2(R'-R)$. If $x,y,a,b$ satisfy
$d(x,a)+d(y,b)\ge
2R-r=\break 2R'-r',\; d(x,y)\le r\le r',\;d(a,b)\le 2N,$ then
$d(y,a)+d(x,b)\le d(x,a)+d(y,b)+2k$. Therefore, $r\in I(R) $.

\vglue4pt c) Let $r\in I(R)$. By b), we may replace $R$ by
$\sup\{R,N\}$ and
$r$ by $\inf \{r,R+N\}$. We assume that $r>R-N+k+6\delta$,
and show that
$r\not\in I(R)$.

Let $y,a'\in X$ be such that
$d(y,a')\ge R+N$. Since $X$ is weakly $\delta$-geodesic, there
exists $a\in X$ such that $d(y,a)\le R+N+2\delta$ and
$d(a,a')\le d(y,a')-\break (R+N)$, whence $R+N\le d(a,y)\le R+N+2\delta$.
Choose $x,b\in X$ such that $d(x,y)\le r\,,$
$d(x,a)\le d(y,a)-r+2\delta$,
$d(a,b)\le 2N$ and $d(y,b)\le d(y,a)-2N\break +2\delta$. We have:
\begin{eqnarray*}
d(x,a)+d(y,b)&\ge &2d(y,a)-d(x,y)-d(a,b)\\
&\ge& 2(R+N)-r-2N \ge 2R-r.
\end{eqnarray*}

Now note that 
$
d(x,b) \ge  d(y,a)-d(y,b)-d(x,a)$, $d(y,a)-d(x,a) 
 \ge 
r-2\delta$ and $d(y,a)-d(y,b)\ge 2N-2\delta$; whence
\begin{eqnarray*}
d(y,a)+d(x,b)-d(x,a)-d(y,b)&\ge& 2(d(y,a)-d(x,a)-d(y,b))\\
&\ge&
2(r+2N-4\delta-d(y,a))
\\
&\ge &2(r+2N-4\delta -R-N-2\delta)>2k.
\end{eqnarray*}
 Therefore
$r\not\in I(R)$.
\enddemo

 For $x\in X$ and $S\in \Delta$, put $A_{S,x}=\{\, a\in
U_S\,,\;d(x,a) \le d(x,U_S)+2\delta\,\}$ and $C_{S,x} = \{\, c\in S\,,
\; d(x,c) \ge \max\{\, d(x,y)\,,\;y\in S \,\}-2\delta\, \}$.

For $x\in X\,,\; S\in \Delta$ and $r\in{\bf R}_+$, we set
$$Y_{S,x,r} = \dst\bigcup_{y\in \overline B(x,r)} A_{S,y} \quad \hbox{and}
\quad Z_{S,x,r} = \dst\bigcup_{y\in \overline B(x,r)} C_{S,y}.$$

\vglue12pt {\it Remark}.  In order to construct the element
$\gamma$, we just need to consider the sets $A_{S,x}$ and
$Y_{S,x,r}$; the sets
$C_{S,x}$ and $Z_{S,x,r}$ are used in the construction of
the ``dual Dirac element'' which will be given in the next section.

\proclaim{Lemma} Assume that $r\in I(d(x,U_S))$.
\begin{itemize}
\ritem{a)} For all $a\in
Y_{S,x,r}${\rm ,} $d(a,x)\le d(x,U_S)+2k+2\delta $. For all $b\in
Z_{S,x,r}${\rm ,}  $d(b,x)\ge \max\{\, d(x,y)\,,\;y\in S
\,\}-2k-2\delta $.
\ritem{b)} The diameter of $Y_{S,x,r}$ is $\le N$. 
 \ritem{c)} For all $y\in \overline B(x,r)$ and all $a\in Y_{S,x,r}${\rm ,} 
$d(a,y)\le d(y,U_S)+4k+2\delta $. 
 \ritem{d)} If moreover $0 <\sup\{\,d(x,c)\,,\; c\in S\,\}-4k-6\delta ${\rm ,} then
$Y_{S,x,r}\cap Z_{S,x,r}=\emptyset$ and the distance between
$Y_{S,x,r}$ and $Z_{S,x,r}$ is $\ge 2\delta $.
\end{itemize}

\endproclaim

\demo{Proof} a) Choose $y,z\in \overline B(x,r)$ such that $a\in A_{S,y}$ and
$b\in C_{S,z}$. Let also $a'\in U_S\,,\;b'\in S$ be such that
$d(x,a')= d(x,U_S)$ and $d(x,b')=\max\{\, d(x,c)$, $c\in S \,\}$.

Note that $d(x,y)\le r\,,\; d(y,a)\ge d(x,U_S)-r\,,\; d(x,a')=d(x,U_S)$,
the diameter of $U_S$ is $\le 2N$, and $ d(y,a)\le
d(y,a')+2\delta$. As $r\in I(d(x,U_S))$ we have $d(x,a)\le
d(x,a')+d(y,a)-
d(y,a')+2k\le d(x,a')+2k+2\delta$.

Also $d(x,z)\le r\,,\; d(y,b')\ge d(x,U_S)-r\,,\; d(x,b)\ge d(x,U_S)$, the
diameter of $S$ is $\le N$, and $ d(z,b')\le d(z,b)+2\delta$. As
$r\in I(d(x,U_S))$ we have $d(x,b')\le d(x,b)+2k+d(z,b')- d(z,b)\le
d(x,b)+2k+2\delta$.

\vglue4pt b) By a), $Y_{S,x,r}\subset \{\, a\in U_S\,,\;d(x,a) \le
d(x,U_S)+2k+2\delta\,\}$. Hence, b) follows from Lemma~6.2.

\vglue4pt c) Choose $a''\in U_S$ such that $d(y,a'')=d(y,U_S)$. By
a), we have $d(x,a)\le d(x,U_S)+2k+2\delta\le d(x,a'')+2k+2\delta$. As
$r\in I(d(x,U_S))\,,$ it follows that $d(y,a)+d(x,a'')\le
d(y,a'')+d(x,a)+2k $, whence
$d(y,a)\le d(y,U_S)+4k+2\delta $.

\vglue4pt d) follows from a) and Lemma~6.4.
\enddemo

\proclaim{Proposition} Let $x\in X${\rm ,} $S\in \Delta$ and $r\in I(d(x,U_S))$.
Then
\begin{itemize}
\ritem{a)} $S \cup Y_{S,x,r}\in \Delta$. 
\end{itemize}
\indent Let $T\in \Delta$ be such that $S-Y_{S,x,r} \subset T
\subset S\cup Y_{S,x,r}$. Then\/{\rm :}\/ 
\begin{itemize}
\ritem{b)} For all $y\in
\overline B(x,r)${\rm ,} $A_{S,y} = A_{T,y}${\rm ,} and
$Y_{S,x,r} = Y_{T,x,r}$.
\ritem{c)} If moreover
$r<\sup\{\,d(x,c)\,,\; c\in S\,\}-4k-6\delta ${\rm ,} then for all $y\in
\overline B(x,r)${\rm ,} $C_{S,y} = C_{T,y}$, and
$Z_{S,x,r} = Z_{T,x,r}$. \end{itemize}

\endproclaim

\demo{Proof} a) is a consequence of Lemma~6.6.b).

\vglue4pt b) It follows from~6.6.c) that for every $a$ in the
symmetric difference of $S$ and $T$ and every $y\in \overline
B(x,r)$ we have: $d(y,a)\le d(y,U_S)+4k+2\delta$. It remains to apply~6.3.b).

\vglue4pt In the hypothesis of c), $\sup\{\,d(y,c)\,,\;c\in
S\,\}>4k+6\delta$; by Lemma~6.4, each point of $C_{S,y}$ is at a
distance $>d(y,U_S)+4k+4\delta$ from $y$. Therefore,
$C_{S,y}\subset T$ by~6.3.b). In the same way, $C_{T,y}\subset
S$; whence
$C_{S,y} = C_{T,y}$.
\enddemo

 \vglue4pt
\centerline{\bf Construction of $\gamma$}
\vglue6pt
Let $\delta ,k,N$ be positive real numbers such that
$N>8k+22\delta $ and let $(X,d)$ be a weakly
$\delta$-geodesic, locally finite metric space satisfying conditions
(C2) and the following:

\begin{itemize}
\item[(C1)] $\;\dst\bigcup_{R\in\R_+} I(R) = \R_+$.
\end{itemize}
Note that, if $\delta\le k$, condition (C1) is slightly weaker than
condition (B1).

Assume now that $X$ has bounded coarse geometry, more precisely,
that:
 \begin{itemize}
\item[(C3)] There exists a $\Gamma $-invariant measure $\mu $
on
$X$ with the property that for any $R>0$ there exists $K>0$ such that for any $x\in X\,,\;\mu
(\overline B(x,R))\le K$ and $\mu (B(x,\delta ))\ge 1$.
\end{itemize}

For a nonempty subset $T$ of $X$, let $\mu _T$ be the
measure defined in the following way: let $\widetilde T$ be the $\delta
$-neighborhood of the set $T$ in $X$; for $x\in \widetilde T$ let
$p_x$ be the probability measure equi-distributed on all points of
$T$ minimizing the distance from $x$ to $T$. Set then $\mu
_T=\int_{\widetilde T} p_x\,d\mu (x)$.

The map $T\mapsto \mu _T$ is $\Gamma $-equivariant and
satisfies:

\proclaim{Lemma} {\rm a)} For every $x\in T${\rm ,} $\mu _T(B(x,2\delta)\cap T)\ge
1$. 
\vglue4pt
 {\rm b)} Let $T$ and $T'$ be subsets of $X$. The measures $\mu
_T$ and $\mu _{T'}$ coincide outside the $2\delta$\/{\rm -}\/neighborhood
of the symmetric difference of $T$ and $T'$.
\endproclaim

\demo{Proof} For every $x\in T$ and any $y\in B(x,\delta)$, the measure
$p_y$ is supported in $B(x,2\delta)$. Since $B(x,\delta)\subset
\widetilde T$ and $\mu (B(x,\delta))\ge 1$, we deduce a).

\vglue4pt b) For every $y\in T$, the value $\mu_T(y)$ only
depends on the measures $p_x$ for $x\in B(y,\delta)$. Moreover,
for $x\in
\widetilde T$, the measure $p_x$ is nonzero only on $T\cap
B(x,\delta )$.
\enddemo

For $x\in X\,,$ and $ S\in \Delta$ put $r_{S,x} =
\sup I(d(x,U_S))-k$.

\proclaim{Proposition} For any $x\in X$ and $S\in \Delta${\rm ,} there exist probability
measures $\psi_{S,x}$ and $\theta_{S,x}$ on $X$
which depend $\Gamma $\/{\rm -}\/equivariantly on the pair $(S,x)$
and have supports in $\dst\bigcup_{r\in I(d(x,U_S))}Y_{S,x,r-k}$ and
$\dst\bigcup_{r\in I(d(x,U_S))} Z_{S,x,r-k}$ respectively such that\/{\rm :}
\begin{itemize}
\ritem{a)} For any $x,y\in X,$ the functions $S\mapsto
\|\psi_{S,x}-\psi_{S,y}\|_1$ and $S\mapsto
\|\theta_{S,x}-\theta_{S,y}\|_1$ converge to zero outside finite sets
of $ \Delta$.
\ritem{b)} For any $x,y\in X$ such that
$d(x,y)\le k${\rm ,} and any $S,T\in
\Delta$ whose symmetric difference is contained in the support of
$\psi_{S,x}${\rm ,} we have $\psi_{S,y}=\psi_{T,y}$\/{\rm ;} if moreover
$d(x,U_S)\ge 6\delta${\rm ,} then $\theta _{S,x}=\theta_{T,x}${\rm ,} and the
distance between ${\rm supp}\, \psi_{S,x}$ and ${\rm supp}\, \theta_{S,x}$ is $\ge
2\delta$.
\end{itemize}

\endproclaim

\demo{Proof} Let $f_{S,x,r}$ (resp.\ $g_{S,x,r}$) denote the characteristic
function of $ Y_{S,x,r}$ (resp.\ $Z_{S,x,r}$). Put 
\begin{eqnarray*}
\widetilde \psi
_{S,x}&=& (f_{S,x,0}+\dst\int_{0}^{r_{S,x}}\!\!f_{S,x,t} \,\,dt\;)\cdot
\mu_{U_S}\;,\\
\widetilde \theta _{S,x}&= &(g_{S,x,0}
+\dst\int_{0}^{r_{S,x}}\!\!g_{S,x,t}\,\,dt\;)\cdot
\mu_{S}\,.
\end{eqnarray*}
Define $\psi_{S,x}$ and $\theta_{S,x}$ to be
proportional to $\widetilde \psi_{S,x}$ and
$\widetilde\theta_{S,x}$ and
normalized to mass $1$.

In order to prove property a), we will first prove that (for fixed
$x,y$), the functions $S\mapsto\|\widetilde \psi_{S,x}-\widetilde
\psi_{S,y}\|_1$ and $S\mapsto\|\widetilde \theta_{S,x}-\widetilde
\theta_{S,y}\|_1$ are bounded.

Set $d=d(x,y)$ and let $K$ denote the maximal $\mu $-volume
of balls of radius $N+\delta$ in $X$. For every subset $Y$ of
$X$, the $\mu_Y$-volume of any ball of radius $N$ in $Y$ is
$\le K$.

Note that $r \mapsto Y_{S,x,r}$ and $r \mapsto Z_{S,x,r}$ are
increasing
functions; therefore, $r \mapsto f_{S,x,r}$ and $r \mapsto g_{S,x,r}$
are also
increasing functions. Their norm is bounded by $K$ because $S$ and
$U_S$
sit inside a ball of radius $N$ with center at any point of $S$.
Also note that $\overline B(y,r) \subset \overline B(x,r+d)$. So by
definition of
$Y_{S,x,r}$ and $Z_{S,x,r}$ we have $ Y_{S,y,r}
\subset Y_{S,x,r+d}$ and $Z_{S,y,r} \subset Z_{S,x,r+d}$. Therefore,
$f_{S,y,r} \le f_{S,x,r+d}$ and $g_{S,y,r}
\le g_{S,x,r+d}$. Using the fact that $r_{S,y}\le r_{S,x}+2d$ (Lemma~6.5.b), we find $\widetilde \psi_{S,y}\le \widetilde
\psi_{S,x}+f_{S,y,0}+3df_{S,y,r_{S,y}}$; exchanging the roles of $x$
and $y$, we deduce
$\|\widetilde \psi_{S,x}-\widetilde
\psi_{S,y}\|_1\le (3d+1)K\,;$ in the same way, $\|\widetilde
\theta_{S,x}-\widetilde
\theta_{S,y}\|_1\le (3d+1)K\,$.

Let $a\in U_S$ and $b\in S$ be such that $d(x,a)=d(x,U_S)$ and
$d(x,b)=\sup\{\, d(x,y)\,,\allowbreak \;y\in S\,\}$. Note that the ball
in $U_S$
(resp.\ $S\,)$ with center $a$ (resp.\ $b$) and radius $2\delta$
is contained in $A_{S,x}$ (resp.\ $C_{S,x}$). By Lemma~6.8.a),
$\mu_{U_S}(A_{S,x})\ge 1$ and $\mu_S(C_{S,x})\ge 1$. Therefore
the norms of $\widetilde \psi_{S,x}$ and $\widetilde \theta_{S,x}$ are
$\ge r_{S,x}$. Now, condition (C1) implies $\lim_{S\ra \infty}
r_{S,x}=\infty$; property a) follows.

We now come to property b). Choose among $ \{x,y\}$ the point which
has the maximal distance to $U_S$ and call it $z$. As the function
$R\mt I(R)$ is increasing and as
$d(x,z)\le k$, every point in the support of $\psi_{S,x}$ is in some
$Y_{S,x,r}\subset Y_{S,z,r+k}$ where $r+k\in I(d(x,U_S))\subset
I(d(z,U_S))$. Now $y\in \overline B(z,k)$; moreover, for any $y'\in
X$
such that $y=y'$ or $d(y,y')<r_{S,x}$, we have $$d(y',z)\le
k+d(y,y')\in I(d(y,U_S))\subset I(d(z,U_S)).$$ By Proposition~6.7.b),
$A_{y',S}=A_{y',T}$; also, the measures $\mu_{U_S}$ and
$\mu_{U_T}$ coincide in $\{\,c\in U_S\,,\; d(z,c)\le
d(z,U_S)+4k+4\delta\}$ by Lemmas~6.6.a),~6.3.b) and~6.8.b); we
deduce that $\psi_{S,y}=\psi_{T,y}$.

If $d(x,U_S)\ge 6\delta$, then, by Lemma~6.5.c), for all $r\in
I(d(x,U_S))$ we have $r\le d(x,U_S)+k$; therefore, $r_{S,x}\le
d(x,U_S)$. Note that since $x\not\in U_S$, we have
$\sup\{\,d(x,c)\,,\;\allowbreak c\in S\,\}>N>4k+6\delta$; by Lemma~6.4, we
deduce that $r_{S,x}<\sup\{\,d(x,c)\,, \;c\in S\,\}  -4k-6\delta$;
therefore, by Proposition~6.7.c) the functions
$g_{S,x,r}$ and $g_{T,x,r}$ coincide for
$r<r_{S,x}$. Now, by Lemmas~6.6.d) and~6.8.b), we deduce that
$\theta _{S,x}=\theta_{T,x}$. The last assertion of the proposition
also follows.
\enddemo

Let $\lambda \in \gtM_k$. Set $\psi_{S,\lambda } =
\dst\int_{X}\psi_{S,x}\,\,d\lambda (x)
$. Finally (identifying measures on $X$ and elements in $\ell^1(X)$ --
since $X$
is discrete), let
$\phi_{S,\lambda }\in
\ell^2(X)$ be the function defined by $\phi_{S,\lambda }(y) =
(\psi_{S,\lambda }(y))^{1/2}$.

Denote by $(e_x)_{x\in X}$ the canonical basis of $\ell^2(X)$. Let
$H$ be the Hilbert subspace of $\Lambda^*(\ell^2(X))$ spanned by
$e_{x_1}\wedge\ldots \wedge e_{x_n}$ where $\{x_1,\ldots
,x_n\}$ runs over $\Delta$. Let $e_\emptyset $ be the vacuum
vector of $\Lambda^*(\ell^2(X))$. The grading of $H$ is the
opposite to \pagebreak the canonical one; \ie the degree of $e_{x_1} \wedge
\ldots \wedge e_{x_n}$ is $n-1$ (modulo~$2$). For $\xi
\in\ell^2(X)$ let $c(\xi)$ be the operator on $\Lambda^*
(\ell^2(X))$ given by the Clifford multiplication by $\xi$, \ie 
$c(\xi)={\rm ext}(\xi)+{\rm int}(\xi)$, the sum of operators of (left) exterior and 
interior multiplication by $\xi$.

Denote by $P$ the projection of $\Lambda^*(\ell^2(X))$ onto
$H$ and let $F_\lambda $ be the operator on $H$ given for
$\{x_1,\ldots,x_n\} = S\in \Delta$ by the formula:
\smallbreak\centerline{${\displaystyle
 F_\lambda (e_{x_1} \wedge \ldots \wedge e_{x_n}) = P\,
c(\phi_{S,\lambda })(e_{x_1} \wedge \ldots \wedge e_{x_n})\,.}$}

\proclaim{Theorem} {\rm a)} For all $\lambda \in \gtM_k\,,\; F_\lambda = F_\lambda ^*$
and $1-F_\lambda ^2$ is the rank\/{\rm -}\/one projection onto the vector
$\phi_{{\rm support}(\lambda),\lambda}$.
 
{\rm b)} For all $\lambda
,\lambda '\in \gtM_k\,,\; F_\lambda - F_{\lambda '}
\in {\cal K}(H)$. 

\vglue4pt
{\rm c)} The map $\lambda \mapsto
F_\lambda $ is norm\/{\rm -}\/continuous on
$\gtM_k$.
\endproclaim

{\it Proof}. If $S=\{x_1,\ldots x_n\}$ is a finite subset of $X$ denote by
$L_S$ the line in $\Lambda^*(\ell^2(X))$ spanned
by $e_{x_1}\wedge
\ldots\wedge e_{x_n}$.
\vglue4pt
  a) Let $\lambda \in \gtM_k$ and $S\in \Delta$. Let
$x$ be a point in the support of $\lambda $ maximizing the
distance to $U_S$. Since the function $R\mapsto I(R)$ is increasing,
the support of
$\phi_{S,\lambda }$ is contained in
$\dst\bigcup_{r\in I(d(x,U_S))}Y_{S,x,r}$.
\vglue2pt
For all $S\in \Delta,\;c(\phi_{S,\lambda })L_S$ is contained in the
sum of spaces $L_T$ where $T$ runs over all subsets of $X$
such that the symmetric difference of $S$ and $T$ consists of
exactly one point which is contained in the support of $\phi_{S,\lambda }$;
in particular, $T\in \Delta \cup\{\emptyset\}$, and for every $y$ in
the support of $\lambda
$ we have $\psi_{S,y} =\psi_{T,y}$ if $T\ne \emptyset$ (Proposition~6.9.b). Therefore
$c(\phi_{S,\lambda })(L_S)\subset H\oplus {\bf C}e_\emptyset$
where $e_\emptyset $ is the vacuum vector of $\Lambda^*( \ell^2
(X))$, and since $c(\phi _{S,\lambda })$ is self-adjoint, we deduce that
$F_\lambda =F_\lambda ^*$.

Consider the equivalence relation in $\Delta $ for which $S$ and
$T$ are equivalent if their symmetric difference is contained in the
support of $\phi_{S,\lambda }$. The vector
$\phi_{S,\lambda }$ is constant on the equivalence classes. The
Hilbert space $H$ breaks into an orthogonal sum of finite-dimensional
subspaces spanned by the lines $L_T$ of equivalence classes, $F_\lambda
$ preserves this decomposition and coincides with $c(\phi _{S,\lambda
})$ on any such class except for the class formed by all subsets $S\subset
X$ such that $S\subset {\rm supp}(\phi_{S,\lambda})$.

The lines $L_S$ for all $S$ in this class span a finite-dimensional subspace
$V$ of $H$. The operator $c(\phi_{S,\lambda })$ maps $V\oplus {\bf
C}e_\emptyset$ isomorphically onto itself, so being compressed to $V$ this
operator has a one-dimensional kernel. In fact, $V\oplus {\bf
C}e_\emptyset$ is isomorphic to a finite-dimensional exterior algebra and
$c(\phi_{S,\lambda })$ is a Clifford multiplication operator on it. The
operator $1-F_\lambda^2$ is $0$ on the orthogonal complement of $V$ in $H$,
and on $V$ it is a one-dimensional projection onto the image of ${\bf
C}e_\emptyset$ under the Clifford multiplication by $c(\phi_{S,\lambda })$,
\ie it is the one-dimensional projection onto \pagebreak the vector $\phi_{S,\lambda }$.

It remains to note that all nonempty subsets of the support of
$\phi_{S,\lambda}$ belong to our distinguished class for any $S\subset
{\rm supp}(\lambda)$. This means that the function $\phi_{S,\lambda }$ for this
class is equal to $\phi_{{\rm support}(\lambda),\lambda}$, which establishes
the expression for $1-F_\lambda ^2$.

\vglue4pt b) and c) Let
$\lambda,\lambda'\in \gtM_k$; for $S\in
\Delta$, put $\phi_S=\inf\{\phi_{S,\lambda },\phi_{S,\lambda '}\}$,
and let $F$ be given for $\xi\in L_S$ by the formula $F(\xi) = P\,
c(\phi_S)(\xi)\,.$ Now, $F_\lambda-F$ and
$F_{\lambda'}-F$ are block diagonal. As
$\|\phi_{S,\lambda }-\phi_S\|^2_2\le\|\phi_{S,\lambda
}-\phi_{S,\lambda '}\|^2_2 \le \|\psi_{S,\lambda }-\psi_{S,\lambda
'}\|_1$, it follows from Proposition~6.9.a) that $F_\lambda -F$ is
compact; in the same way,
$F_{\lambda '}-F$ is compact. Moreover,
$\|F_\lambda-F\|^2=\sup\{\|\phi_{S,\lambda
}-\phi_S\|^2_2,\;S\in\Delta \}\le \sup\{\|\psi_{S,\lambda }-\psi_{S,\lambda
'}\|_1,\;S\in\Delta \}\le \|\lambda-\lambda '\|_1$; estimating (in the same
way)
$\|F_{\lambda '}-F\|^2$, we obtain c).
\hfill\qed\vglue8pt

The pair $(H,F_x)$ defines an element $\gamma =
\gamma_k$ of $R(\Gamma )$.

\proclaim{{C}orollary}  
Let $q$ be the map $\gtM_k \ra {\rm point}${\rm .}  Then $q^*(\gamma) = 1$ in $RK^0_{\Gamma }(\gtM_k)$.

\demo{Proof} By Theorem~6.10.a) and c), the family $\lambda \mapsto
F_\lambda $ defines the unit element of $RK^0_{\Gamma }(\gtM_k)$.
By Theorem~6.10.b), it is equal to
$q^*(\gamma_k)$.
\enddemo

\vglue-20pt
\section{The dual Dirac construction}
\vglue-6pt

In this section and the next one we finish the proof of Theorem 5.2 (and
therefore of Theorem 1.1).

Recall that having a proper isometric action of a discrete group
$\Gamma $ on a weakly bolic, weakly geodesic metric space of bounded
coarse geometry, we can choose a discrete, locally finite subspace $X$
of this space having the same properties and still equipped with a proper
isometric action of $\Gamma $ (Prop.~3.5).

\vglue12pt\centerline{{\bf The dual Dirac construction for euclidean spaces}
([HKT], [HK])}
\vglue6pt
Let us denote the real $L^2(X;\mu)$ by $\H$. Recall first the construction
of a \Cst-algebra associated to a separable real Hilbert space ([HKT],
[HK]). We will consider $\H$ as an affine Hilbert
space. For any finite-dimensional real affine subspace
$V\subset \H$, denote by $V_0$ the underlying vector space. Let $C_\tau
(V)$ be (as in [K2], 4.1) the ${\bf Z}_2$-graded algebra of Clifford
functions
$C_0(V,{\rm Cliff}(V_0))$, where ${\rm Cliff}(V_0)$ is the complex Clifford algebra
associated with the Euclidean quadratic form on $V_0$. (This algebra
has a natural structure of a \Cst-algebra.) If we fix an origin $v_0$ in
$V$ (\ie if we identify $V$ with $V_0$), then we can define the {\it Bott
operator\/}
$B_{v_0}$ for $C_{\tau}(V)$ as multiplication by the function $C: V\to {\rm Cliff}
(V_0)$ whose value at
$v\in V$ is $v-v_0\in Cliff (V_0)$. It is a degree-one, essentially
self-adjoint, unbounded multiplier of $C_{\tau}(V)$, with domain the
compactly supported functions in $C_{\tau}(V)$.

Denote by $\cS=C_0({\bf R})$ the ${\bf Z}_2$-graded algebra of continuous
complex-valued functions on
${\bf R}$ which vanish at infinity. $\cS$ is graded according to even and
odd functions. We will use the following notation:
$\A(V)=\cS  \,\widehat{\otimes}\,C_{\tau}(V)$. Denote by $X$ the operator of
multiplication by $x$ on
${\bf R}$, viewed as a degree one, essentially self-adjoint, unbounded
multiplier of $\cS$ with domain the compactly supported functions in
$\cS$.

Suppose now that $V$ and $V'$ are finite-dimensional affine subspaces
of $H$ with $V'\subset V$. Let $W$ be the orthogonal complement of
$V'$ in $V$. Since $W$ naturally has an origin (call it $w_0$),
we can identify $W$
with $W_0$. Define a homomorphism $\cS\to \A(W)$ by the formula:
$f\mapsto f(X \,\widehat{\otimes}\,1+1 \,\widehat{\otimes}\,B_{w_0})$. Using a canonical
decomposition: $\A(V)=\A(W) \,\widehat{\otimes}\,C_{\tau}(V')$, define a
homomorphism $\A(V')\to \A(V)$ by tensoring the above homomorphism
$\cS\to \A(W)$ with the identity map on $C_{\tau}(V')$. These
homomorphisms $\A(V')\to \A(V)$ are transitive with respect to the
embeddings $V''\subset V'\subset V$. The algebra $\A(\H)$ is defined as
the inductive limit of the algebras $\A(V)$ taken over the directed set
of all finite-dimensional affine subspaces of $\H$. We will consider this
\Cst-algebra as trivially graded.

Note that for any point $v\in \H$, there exists a canonical unbounded
multiplier of $\A(\H)$, the unbounded Bott element $B_v$. It corresponds to
the unbounded multiplier $X$ of the algebra $\A(v)=\cS$ (where $v$ is
considered as an affine subspace of $\H$) under the inductive limit map
$\A(v)\to \A(\H)$. The unbounded multiplier $B_v$
is essentially self-adjoint; its resolvent belongs to
$\A(\H)$;
moreover
$B_v-B_w$ is bounded for $v,w\in \H$, and $gB_v=B_{gv}$ for all
$g\in \Gamma$. Therefore $B_v$ defines an element $\eta_\H\in
KK_1^\Gamma (\C,\A(\H))$. Note that in our case of $\H=L^2(X;\mu)$,
the action of $\Gamma$ on $\H$ is
{\it linear}, so the operator $B_0$ (corresponding to $v=0$) is exactly
$\Gamma$-invariant.

\vglue16pt \centerline{\bf The algebra $\A_k$}
\vglue6pt

In the previous section, for any $k\in \R_+$, we have defined $N\in \R_+$.
Let us fix $k$ and $N$ and consider the map $\tau :\gtM_N \to L^2(X;\mu)$
defined in 4.2. Denote the image of this map (with the topology induced
by the weak topology of $L^2(X;\mu)$) by $Z_k$. According to Lemma 4.3,
$Z_k$ is a proper $\Gamma $-space.

We will define now a \Cst-algebra $\A_k$ which is a
$\Gamma -Z_k$-algebra. Note that whereas
$\gtM_N$ is obviously a simplicial complex, with simplices
corresponding to subsets $S\subset X$ of diameter $\le N$, the space
$Z_k$ has only a structure of a polyhedral complex. For any simplex $S$
in $\Delta _N$, we denote by $|S|\subset \gtM_N$ its geometric
realization (\ie the set of probability measures on $X$ with support in
$S$); the subspace $\tau(|S|)$ is a convex polyhedron in $L^2(X;\mu)$.
We will denote this polyhedron by $\widetilde S$ and call it {\it
quasisimplex}.

Let $q_{S}$ be the map of
$\H$ into itself which associates to each point of $\H$ the nearest point
of its finite-dimensional subquasisimplex $\widetilde S$. It is clear that
this map is well defined and continuous in the weak topology. Using this
map, one can lift any continuous function on $\widetilde S$ to a weakly
continuous function on $\H$. This map defines a homomorphism of the
algebra of continuous functions on the quasisimplex into the algebra of
weakly continuous functions on $\H$. But since multiplication by a
bounded weakly continuous function on $\H$ naturally defines a central
multiplier of $\A(\H)$ (of degree $0$), we obtain a homomorphism
$C(\widetilde S)\to {\cal M}(\A(\H))$.

Let $H$ be the (graded) Hilbert space entering the definition of the
$\gamma$-element of the previous section.  The \Cst-algebra
$\K(H)\otimes \A(\H)$ is endowed with the natural (diagonal)
action of the group $\Gamma $. Moreover, there exists also a $\Gamma
$-equivariant homomorphism of the algebra of functions
$C_0(Z_k)$ to ${\cal M}(\K(H)\otimes \A(\H))$. Namely, if we fix a
natural basis in $H$
as  in the previous section (each basis vector corresponding
to an oriented simplex of $\gtM_N$), then the Hilbert module
$H\otimes \A(\H)$ over $\A(\H)$ becomes a direct sum over
the simplices $S$ of $\gtM_N$ of the algebras $\A(\H)$. For each
simplex $S$, there is a homomorphism of $C(\widetilde S)\to {\cal
M}(\A(\H))$ as defined above. This homomorphism clearly does not depend on
the orientation of $S$, and all these homomorphisms combine together to
give the homomorphism $C_0(Z_k)\to \oplus_S C(\widetilde S)\to {\cal
L}(H\otimes \A(\H))$.

\numbereddemo{Definition}  Let $\A_k$ be the subalgebra of all elements of the\break \Cst-algebra
$\K(H)\otimes \A(\H)$ which commute with the image of the
above homomorphism $C_0(Z_k)\to {\cal M}(\K(H)\otimes \A(\H))$.
There exists a natural restriction of the homomorphism $C_0(Z_k)\to {\cal
M}(\K(H)\otimes \A(\H))$ to a homomorphism $C_0(Z_k)\to {\cal
M}(\A_k)$. This latter homomorphism defines a structure of a
$C_0(Z_k)$-algebra on $\A_k$.
\enddemo

This definition requires some explanation. First of all, the algebra $\A_k$
is very similar to the algebra of a finite-dimensional simplicial complex 
defined in [KS1] and called $\A_X$ there. We could not use that definition
here because in general the Rips complex $\gtM_N$ is not finite-dimensional. 
The main difference between the algebra $\A_X$ of [KS1] and the algebra 
$\A_k$ defined here is that in [KS1] we used the function called {\it type}
to map simplices of our simplicial complex into a Euclidean space, 
whereas here we just have an embedding of   our polyhedral complex 
$Z_k$ into the  Hilbert  space $\H$. 

To make this parallel between the definition of [KS1] and the definition 
given here completely precise, let us consider the Hilbert module
$H\otimes \A(\H)$ over $\A(\H)$ as a direct sum over
the simplices $S$ of $\gtM_N$ of the algebras $\A(\H)$.
The algebra $\A_k$ clearly contains all elements of $\K(H\otimes \A(\H))$
which are diagonal with respect to this decomposition. It also contains 
nondiagonal elements, namely for those pairs of copies of $\A(\H)$ in the
direct sum $\mathbold{\oplus}_S \A(\H)$ which correspond to quasisimplices $\widetilde S$, 
$\widetilde S'$ with $\widetilde S \cap \widetilde S'\ne \emptyset$. One easily 
verifies that $\M(\A_k)$ is a subalgebra of $\L(H\otimes \A(\H))$ and that 
$\A_k=\M(\A_k) \cap \K(H\otimes \A(\H))$.

The algebra $\A_k$  just constructed is a
proper $\Gamma $-algebra (Def.\ 5.6). This will be the algebra $\A_k$ of
Corollary 5.14. In order to apply Corollary 5.14, it is enough to construct
the Dirac and  the dual Dirac elements $d_k$ and
$\eta_k$. The Dirac element will be constructed in the next section. Here we
construct the dual Dirac element $\eta_k\in
K_1^{\Gamma }(\A_k)$. First we will choose for any simplex $S$ of the
Rips complex a certain vector $v_S\in \H$ which will be the center of the
Bott element for this simplex.

\vglue18pt

\centerline{\bf The choice of the vector $v_S$}
\vglue6pt

In the previous section, we have associated to any simplex
$S$ of the combinatorial Rips complex $\Delta_N$ and any point $x\in
X$, a measure $\theta_{S,x}\in \gtM_N$ (Prop.\ 6.9). Let us put
$\xi_{S,x}=\tau(\theta_{S,x})$ where $\tau$ is the map defined in 4.2.

For all $R>0\,,$ let $ K(R)>0$ denote the maximal
$\mu$-measure of a subset $T\subset X$ of diameter $\le R$.
Put $K=K(N+2\delta)$.

\proclaim{Lemma} 
Let $\rho,\alpha$ be positive real numbers such that
$\alpha>2\rho K^{3/2}+K^2$. Put $v_S=\alpha\xi_{S,x}${\rm ,} and let $z$ be any
point of ${\cal H}$ contained in a ball of radius $\rho$ around $v_S$. Suppose
that a point $a\in S$ is at a distance $\ge 2\delta$ in $X$ from the support of
the function $\theta_{S,x}$. Then $q_S(z)$ {\rm (}\/the nearest\/{\rm -}\/to\/{\rm -}\/$z$ point of the
quasisimplex $\widetilde S${\rm )} actually belongs to the subquasisimplex
$(S-\{a\})^{\widetilde{\;\;}}$ in ${\cal H}$.
\endproclaim

\demo{Proof} Put $y=q_S(z)$. The condition that $y$ is the nearest-to-$z$
point of
the quasisimplex $\widetilde S\subset{\cal H}$ means that $\widetilde S$ lies
in the half-space which is cut off from $z$ by a hyperplane passing
through the point $y$ and orthogonal to the vector $y-z$. This can be
expressed by saying that for any point $s\in\widetilde S$, the scalar
product $\langle s-y,y-z\rangle $ is nonnegative. In particular
$\langle \xi_{S,x}-y,y-z\rangle \ge 0$, since $\xi_{S,x}\in
\widetilde S$.

Let $s_1,\ldots ,s_m$ be all vertices of $S$. Then $y$ is a convex linear
combination of the points $\tau(s_1),\ldots ,\tau(s_m)$. The nonzero coefficients
in this linear combination can correspond only to those points among
$\tau(s_1),\ldots ,\tau(s_m)$ that lie in the hyperplane passing through the point
$y$ and orthogonal to the vector $y-z$. We claim that $\langle
\tau(a)-\xi_{S,x},y-z\rangle >0$. It follows that
$\langle \tau(a)-y,y-z\rangle >0$. Therefore $\tau(a)$ cannot lie in this
hyperplane, and the lemma follows.

To prove our claim, note that by Lemma 4.3, $\|\xi_{S,x}\|\ge K^{-1/2}$ and
$\|y\|\le K^{1/2}$. Therefore 
\begin{eqnarray*}
\langle \xi_{S,x},z-y\rangle &=&\langle
\xi_{S,x},v_S+(z-v_S)-y\rangle\ge \alpha
\|\xi_{S,x}\|^2-\|\xi_{S,x}\|(\rho +K^{1/2})\\ &
\ge& \alpha K^{-1}-K^{1/2}(\rho +K^{1/2}).
\end{eqnarray*}

Note also that an element of $\widetilde S$ is a nonnegative function on
$X$, hence $\langle \tau (a),y\rangle \ge 0$.

Since the point $a\in S$ is at a distance $\ge 2\delta$ in $X$ from the
support of the function $\theta_{S,x}$, one has $\langle \tau
(a),v_S\rangle =0$. So 
\begin{eqnarray*}
\langle \tau (a)-\xi_{S,x},y-z\rangle
&=&\langle \tau (a),y\rangle -\langle\tau (a),z-v_S\rangle +\langle
\xi_{S,x},z-y\rangle \\
&\ge& -\rho K^{1/2}+\alpha K^{-1}-K^{1/2}(\rho+
K^{1/2})>0
\end{eqnarray*}
 as claimed.
\enddemo

We choose $\rho=1$ and put $\alpha=2(K^{3/2}+K^{2})$; hence,
$v_S=2(K^{3/2}+K^{2})\xi_{S,x}$.

\vglue18pt

\centerline{\bf The operator $\Phi$}
\vglue6pt

We are ready now to construct an operator $\Phi\in {\cal M}(\A_k)$
which will give us the required dual Dirac element in
$K_1^{\Gamma }(\A_k)$. Let us denote the Hilbert module
$H\otimes \A(\H) $ over $\A(\H)$ by
${\frak H}$. We will consider ${\frak H}$ as a direct sum over the
oriented simplices
$\{S\}$ of $\gtM_N$ of the algebras $\A(\H)$. The operator $\Phi$ will
be constructed in
$\L({\frak H})={\cal M}(\K(H)\otimes \A(\H))$ and then we will prove
that it actually belongs to ${\cal M}(\A_k)$.

Let $B_{v_S}$ be the unbounded Bott element for $\A(\H)$ centered at the
point $v_S$. Denote by $\dim [S]$ the dimension of the simplex
$S$ in $\gtM_N$. Define a function $\nu$ on the real line by
$\nu(t)=t\cdot (\max(1,\vert t\vert))^{-1}$. Also let $\omega$ be a
positive continuous function on $\R_+$ satisfying the conditions:
$\omega (t)=0$ if $t\le 6\delta$ and $\omega(t) =1$ if $t\ge 8\delta$.

The operator $\Phi$ consists of a part $\Phi_{\rm diag}$ diagonal with
respect to the above direct sum decomposition of ${\frak H}$ and an
off-diagonal part $\Phi_{\rm off}$. The operator $\Phi$ will depend on the
point $x\in X$. When we need to show the dependence of
$\Phi$ on $x$ we will denote it by $\Phi_{x}$.

Recall from the proof of Theorem 6.10 that if $S=\{x_1,\ldots x_n\}$ is
a finite subset of $X$, we denote by
$L_S$ the line in $\Lambda^*(\ell^2(X))$ spanned
by$e_{x_1}\wedge\ldots\wedge e_{x_n}$.
The diagonal part $\Phi_{\rm diag}$ of our operator is defined, for
$S=\{x_1,\ldots ,x_n\}\in
\Delta_N$ and $\xi\in L_S$,
by $\Phi_{\rm diag}(\xi\otimes \zeta)= (-1)^{\dim [S]-1}\xi\otimes
\nu (B_{v_S})(\zeta)$.

The off-diagonal part $\Phi_{\rm off}$ takes into account the operator $F=F_x$
of the $\gamma$-element constructed in the previous section (here $x\in X$
is considered as the Dirac measure $\lambda$ at the point $x$). Let
$\widetilde\omega$ be the multiplication by $\omega(d(x,U_S))$ on
$L_S\otimes \A(\H)$. We set  $\Phi_{\rm off}=(F\otimes
1)(1-\Phi_{\rm diag}^2)^{1/2}\widetilde\omega $, and 
finally,    $\Phi=\Phi_{\rm diag}+\Phi_{\rm off}$.

\proclaim{Theorem}  The element $\Phi$ is a multiplier of the algebra $\A_k$. It satisfies
the following conditions\/{\rm :} $\Phi^*=\Phi,\;1-\Phi^2 \in \A_k$. The element
$\Phi$ depends equivariantly on the point $x \in X${\rm ,} namely\/{\rm :}
$g(\Phi_x)=\Phi_{gx}${\rm ,} for all $g\in \Gamma ${\rm ,} and for any $y \in X,\quad
\Phi_{y}-\Phi_{x}\in {\cal A}_k$. Thus $\eta_k =({\cal A}_k,\Phi)$ is
an element of $K_1^{\Gamma }(\A_k)$. Moreover{\rm ,} the image of $\eta_k$ in
the group $K_1^\Gamma(\K(H)\otimes\A(\H))$ under the natural embedding
$\A_k\subset
\K(H)\otimes\A(\H)$ is equal to the product of the $\gamma$\/{\rm -}\/element $\gamma_k$
defined by the operator $F$ and the Bott element $\eta_\H \in
KK_1^\Gamma(\C,\A(\H))$.
\endproclaim

\demo{Proof} Obviously  the operators $\widetilde\omega$ and $\Phi_{\rm 
diag}$ commute. According
to Lemmas 6.3.a) and 6.6.a), if $S,T\in
\Delta_N$ are such that
$F(L_S)$ and $L_T$ are not orthogonal, then $d(x, U_S)=d(x, U_T)$.
Therefore $\widetilde \omega$ and $F\otimes
1$ commute. If moreover $d(x, U_S)\ge 6\delta$, then according to
Proposition 6.9.b) the vectors
$\xi_{S,x}$ and $\xi_{T,x}$ coincide, so that $v_S=v_T$. It follows that 
$\Phi_{\rm diag}$ and
$(F\otimes 1)\widetilde\omega$ anticommute. Therefore
$\Phi^*=\Phi$ and $1-\Phi^2=(1-\widetilde\omega^2 (F^2\otimes 
1))(1-\Phi_{\rm diag}^2 )
\in \K(H)\otimes \A(\H)$ (recall that the space $X$ is locally 
finite, and the condition
$d(x, U_S)\ge 8\delta$ is not satisfied only for a finite number of 
simplices $S$).

The assertion concerning the equivariance of $\Phi$ is clear.

By Proposition 6.9 and Lemma 4.3, we have:  $\Phi_{{\rm 
diag},y}-\Phi_{{\rm
diag},x}\in \K(H)\break\otimes \A(\H);$  from Theorem 6.10, it follows 
that $$((F_x-F_y)\otimes
1)(1-\Phi_{{\rm diag},x}^2)^{1/2}\in \K(H)\otimes \A(\H).$$  Therefore 
$\Phi_{y}-\Phi_{x}\in
\K(H)\otimes \A(\H)$.

Now we have to prove that the operator $\Phi$ is a multiplier of $\A_k$.
Indeed, the operator $\Phi_{\rm diag}$ is obviously a
multiplier of $\A_k$, so we need to prove only that
$\Phi_{\rm off}$ is also a multiplier of $\A_k$.
Suppose that $F(L_S)$ and $L_T$ are not orthogonal. Then we have to
show that the operators of multiplication by functions from
$C_0(Z_k)$ are the same on the two copies of $\H$ corresponding to $S$
and $T$ in the regions where
$1-\Phi_{\rm diag}^2\ne 0$. These regions are in fact balls of radius $1$
around points $v_S$ and $v_T$ respectively.

First of all, in the case when $d(x, U_S)\le 6\delta$, the operator
$\Phi_{\rm off}$ is zero by definition. Therefore, we may assume that
$d(x, U_S)\ge 6\delta$. According to Lemma 7.2, applied in the case when
$a$ is any point in the symmetric difference of $S$ and
$T$, the projection maps for the ball of radius $1$ around $v_S$ (or
$v_T$) to $\widetilde S$ and to $\widetilde T$ are identical, and the 
image of these maps belongs to the intersection of $\widetilde S$ and
$\widetilde T$. (Note that the condition that the point
$a$ is at a distance $\ge 2\delta$ in $X$ from the support of the
function $\theta_{S,x}$ is satisfied in view of the last assertion of
Proposition 6.9.) Therefore, the operators of multiplication by functions
lifted from $\widetilde S$ and
$\widetilde T$, which are defined using these projections, are the same
in the two balls that we consider.

The last assertion of the theorem
follows from the existence of an obvious homotopy moving all centers
$\{v_S\}$ of the Bott elements used in the construction of the operator
$\Phi$ to the zero point of $\H$. The function $\omega$ in the process
of this homotopy is replaced by $1$.
\enddemo
\vglue-12pt
\section{End of proof of the main result}
\vglue-4pt
In Section 5 we   reduced the proof of Theorem 5.2 to the
verification of the sufficient conditions given in Corollary 5.14. In this
section we prove that these sufficient conditions \pagebreak  are fulfilled for groups
satisfying the assumptions of Theorems 1.1, 5.2. More precisely, we
will construct a
Dirac element $d_k\in KK^\Gamma (\A_k,\C)$ such that (in the notation of
Corollary 5.14),
$p_{\gtM_k}^*(\eta_k\otimes_{\A_k} d_k)\break =1_{\gtM_k}$.
\vglue2pt
Consider the extension of
\Cst-algebras
$$0\ra \K(\E)\ra\F\ra \A(\H)\ra 0 \eqno (1)$$
constructed in [HK], where $\E$ is a Hilbert $D$-module. Here
$D=C_0((0,1)^2)$ since we view all \Cst-algebras as trivially graded.
Let us tensor this extension
with $\K(H)$:
$$0\ra
\K(H)\otimes \K(\E)\ra\K(H)\otimes\F\ra
\K(H)\otimes\A(\H)\ra 0. \eqno (2)$$   

Since the algebra $\A_k$ is proper and nuclear, we can apply Proposition
5.7 to get a six term exact sequence associated with the exact sequence (2).

\numbereddemo{Definition} The Dirac element $d_k\in
KK_1^\Gamma(\A_k,\K(H)\otimes\K(\E))=KK_1^\Gamma(\A_k,\C)$  is defined to be
$\partial (u_k)$, where
$u_k\in KK^\Gamma (\A_k,\K(H)\otimes\A(\H))$ is associated with the embedding
$\A_k\to
\K(H)\otimes\A(\H)$, and $\partial$ is the connecting map of the exact
sequence (2).
\enddemo

\proclaim{Proposition}
$p_{\gtM_k}^*(\eta_k\otimes_{\A_k} d_k)=1_{\gtM_k}$ {\rm (}\/up to sign\/{\rm ).}
\endproclaim

\demo{Proof} The action of $\Gamma$ on $\gtM_k$ is proper. Applying Remark 5.10, we find
$$p_{\gtM_k}^*(\eta_k\otimes_{\A_k} d_k)=p_{\gtM_k}^*(\eta_k)\otimes_{\A_k}
p_{\gtM_k}^*\partial
(u_k)=\partial (p_{\gtM_k}^*(\eta_k\otimes_{\A_k} u_k))$$
where we have denoted by the same letter $\partial$ the connecting maps
$$RKK^\Gamma(\gtM_k;\A_k,\K(H)\otimes\A(\H)) \to
RKK_1^\Gamma(\gtM_k;\A_k,\K(H)\otimes\K(\E))$$ and
$$RKK_1^\Gamma(\gtM_k;\C,\K(H)\otimes\A(\H))\to
RKK^\Gamma(\gtM_k;\C,\K(H)\otimes\K(\E))$$ associated with the exact sequence
(2).

The last assertion of Theorem 7.3 implies that $\eta_k\otimes
_{\A_k}u_k=\gamma_k\otimes_\C \eta_\H$ where $\eta_\H\in
KK^\Gamma_1(\C,\A(\H))$  is the Bott element for
$\A(\H)$. By Corollary 6.11, we then find
$p_{\gtM_k}^*(\eta_k\otimes_{\A_k} u_k)=p_{\gtM_k}^*(\eta_\H)$.
\vglue4pt
Therefore $p_{\gtM_k}^*(\eta_k\otimes_{\A_k}
d_k)=\partial(p_{\gtM_k}^*(\eta_\H))$ where $\partial$ is the connecting map
for the exact sequence (1). Let us denote $C_0(0,1)$ by $S$ and consider the
extension of $\Gamma$-algebras
$$0\ra \K(\E)\ra\F_S\ra S\ra 0 \eqno (3)$$ which is the pull-back of the
extension (1) along the natural embedding $S\to \A(\H)$. (This embedding $S\to
\A(\H)$ is defined by the embedding of the point $0$ into $\H$ -- see
\S7.) If $a$ is the class in $KK^\Gamma_1(\C,S)$ of
the Bott element of $S$, then, according to Remark 5.10, it will be enough to
prove that $\partial(p_{\gtM_k}^*(a))=1_{\gtM_k}$
where $\partial$ now is the connecting map for the exact sequence (3).

In the proof of Theorem 6.10 in [HK], the extension (3) was presented as a sum
of two extensions: one of them was $$0\to S\otimes C_0(0,1)\to S\otimes C_0[0,1)\to S\to 0 \eqno
(4)$$
and the second was a restriction (to the point $1$ of the half-interval
$(0,1]$)
of a certain extension $$0\to C_0(0,1]\otimes \K(\E')\to \F'_S\to S\to 0 \eqno
(5)$$ where $\E'$ is another Hilbert $D$-module. Note that the connecting map
$\partial$ for a sum of two extensions is the sum of connecting maps. For the
extension (4) it is easy to prove directly that (up to sign) $\partial(a)=1$.
On the other hand, for the extension (5), as well as for its restriction to
the point
$1$, the connecting map is $0$ because the algebra $C_0(0,1]\otimes \K(\E')$ is
contractible.
\enddemo

\numbereddemo{{R}emark} It can be proved that in fact $\eta_k\otimes_{\A_k} d_k=\gamma_k$.
Let us indicate the main steps of this proof.
\enddemo

The Dirac element $d_k\in KK_1^\Gamma(\A_k,\C)$ gives rise to an
(equivariant) asymptotic
morphism of $\A_k$ to compact operators in some Hilbert space. Related to
this asymptotic morphism, there exists (see [HK, Def.~7.7]) an
associated homomorphism $(d_k)_*:K_1^\Gamma(\A_k)\ra K_1^\Gamma(\C)$ such
that $(d_k)_*(\eta_k)=\eta_k\otimes_{\A_k} d_k$ (see Theorem 7.8 of [HK]).

The convenience of viewing the latter product as $(d_k)_*(\eta_k)$ is that
$(d_k)_*(\eta_k)$ is functorial in $d_k$. In fact, there exists another
asymptotic morphism, which we will denote $d\otimes 1$, of $\A(\H)\otimes \K(H)$
to the algebra of compact operators such that the asymptotic morphism
corresponding to $d_k$ is its restriction to $\A_k$. (Because the algebra
$\A(\H)\otimes
\K(H)$ is not proper, we cannot define a $KK$-element corresponding to this
asymptotic morphism.) This allows one to replace $(d_k)_*(\eta_k)$ with
$(d\otimes 1)_*(\eta_\H\otimes_\C \gamma_k)$ (by functoriality and Theorem 7.3).
(Note that $(d\otimes 1)_*$ is defined only on elements of the group
$K_1^\Gamma(\A(\H)\otimes \K(H))$ which have some special form; see [HK,
Def.\ 7.7]. The element $\eta_\H\otimes_\C\gamma_k$ satisfies this
condition.)

Now $d\otimes 1$ is really a product of the Dirac asymptotic morphism $d$ for
$\A(\H)$ (see [HK, Def.\ 6.7]) and the identity map on $\K(H)$. So
$$(d\otimes 1)_*(\eta_\H\otimes_\C \gamma_k)=d_*(\eta_\H)\otimes_\C
\gamma_k.$$  (The
proof of this equality is similar to the proof of  Theorem~7.8 in [HK].) It
remains to note that $d_*(\eta_\H)=1\in KK^\Gamma(\C,\C)$ by Theorem~6.10 of
[HK], so that $(d\otimes 1)_*(\eta_\H\otimes_\C \gamma_k)= \pagebreak \gamma_k$.

 \centerline{\bf Appendix}
\vglue12pt

In this appendix we prove that $KK^\Gamma(A,B)\simeq E^\Gamma(A,B)$ if A is a
proper nuclear algebra.

Let us first briefly recall the definition of the functor
$E^\Gamma(A,B)$ and the homomorphism $KK ^\Gamma(A,B)
\lra E^\Gamma(A,B)$:

\vglue4pt  -- An asymptotic morphism from $A$ to $B$ is
given by a
$\Gamma-C[0,1]$-algebra $D$ (where $\Gamma$ acts
trivially on $[0,1]$) together with an extension of
$\Gamma-C[0,1]$-algebras $0\ra B(0,1]\ra D\ra A\ra 0$, where
$C[0,1]$ acts on $A$ through evaluation at $0$. Two such
extensions $0\ra B(0,1]\ra D_1\ra A\ra 0$ and $0\ra B(0,1]\ra
D_2\ra A\ra 0$, are said to give the same asymptotic morphism, if
there is an isomorphism $f:D_1\ra D_2$ with a commuting
diagram
$$\begin{array}{ccccccc}
 0\ra &B(0,1]& \lra  & D_1&\lra &A&\ra 0\cr
&\vert\vert&&f\downarrow\ \ &&\vert \vert&\cr 0\ra &B(0,1]&\lra
& D_2&\lra &A&\ra 0.
\end{array}
$$ 

\vglue4pt  -- If $x:\quad 0\ra B(0,1]\ra D\ra A\ra
0$ is an asymptotic morphism from
$A$ to $B$ and $f:B\ra B_1$ is an equivariant
$*$-homomorphism, there is an asymptotic morphism $f(x):\; 0\ra
B_1(0,1]\ra D_1\ra A\ra 0$ from $A$ to $B_1$ uniquely
defined up to isomorphism by the existence of a $\Gamma-C[0,1]$-
equivariant homomorphism $g:D\ra D_1$ such that the
diagram$$\begin{array}{ccccccc}
0\ra &B(0,1]& \lra  & D&\lra &A&\ra 0\cr
&f\downarrow&&g\downarrow\ \ &&\vert \vert&\cr 0\ra
&B_1(0,1]&\lra & D_1&\lra &A&\ra 0
\end{array}
$$ 
commutes. A homotopy is an
asymptotic morphism from $A$ to
$B[0,1]$. Two asymptotic morphisms
$0\ra B(0,1]\ra D_0\ra A\ra 0$ and $0\ra B(0,1]\ra D_1\ra
A\ra 0$ from
$A$ to $B$ are said to be homotopic if there exists an
asymptotic morphism from $A$ to $B[0,1]$ whose evaluation
at $0$ is $0\ra B(0,1]\ra D_0\ra A\ra 0$ and whose
evaluation at $1$ is $0\ra B(0,1]\ra D_1\ra A\ra 0$.
\vglue4pt 
-- Denote by $[[A,B]]$ the set of homotopy classes of asymptotic
morphisms from $A$ to $B$. Denote by $[[A,B]]_{\rm cp}$ the
set of homotopy classes of asymptotic morphisms from $A$ to
$B$ which admit an equivariant completely positive lifting
$A\ra D$ of norm $1$.
\vglue4pt 
-- Denote by $\K$ the algebra of compact operators on
$L^2(\Gamma \times \N)$ where $\Gamma$ acts on
$\Gamma$ by left translation and trivially on $\N$. Let also
$S$ denote the $C^*$-algebra $C_0(\R)$ (trivially graded,
with a trivial action of
$\Gamma$). We write $SA$ instead of $S\otimes A$. By definition,
$E^\Gamma(A,B)=[[SA\otimes \K,SB\otimes \K]]$.

-- Also, according to [T2], $KK^\Gamma(A,B)=[[SA\otimes \K,SB\otimes
\K]]_{\rm cp}$.
\vglue6pt 
-- The morphism $\Phi:KK^\Gamma(A,B)\ra E^\Gamma(A,B)$ is defined as the
forgetful map. This morphism is obviously natural in $B$ (and  also in
$A$, but we do not use, nor define the naturality in $A$ here).
\vglue6pt 
-- Tensoring with $S\K$ leads to well
defined maps $[[A,B]]\ra E^\Gamma(A,B)$ and $[[A,B]]_{\rm cp}\ra
KK^\Gamma(A,B)$. Tensoring again with $S\K$ leads to isomorphisms
$E^\Gamma(A,B)\ra E^\Gamma(SA\otimes \K,SB\otimes \K)$ and
$KK^\Gamma(A,B)\ra KK^\Gamma(SA\otimes \K,SB\otimes \K)$. These
constructions commute with the forgetful map $\Phi$.
\vglue6pt 
-- The asymptotic morphism $0\ra A(0,1]\ra A[0,1]\ra A\ra 0$
defines the unit elements $1_A^K$ of $KK^\Gamma(A,A)$ and
$1_A^E$ of $E^\Gamma(A,A)$. Therefore $\Phi (1_A^K)=1_A^E$.
\vglue6pt
If $E$ is a $C_0(X)$-algebra and $Y\subset X$ is a closed
subset, we denote by $E_Y$ the $C_0(Y)$-algebra
$E/C_0(X-Y)E$.}

\specialnumber{A.1} \proclaim{Lemma}
Let $x:  \; 0\ra B(0,1]\ra
D\ra A\ra 0$ be an asymptotic moprphism. Denote by $p_0:D\ra A$ and
$p_1:D\ra B$ the evaluation maps of the $C[0,1]$\/{\rm -}\/algebra $D$. There
exists an asymptotic morphism
$y:  0\ra D(0,1]\ra \overline D\ra A\ra 0$ such that
$p_0(y)$ is the asymptotic morphism $0\ra A(0,1]\ra A[0,1]\ra
A\ra 0$ and $p_1(y)=x$. If moreover
$p_0$ admits an equivariant completely positive lifting of norm $1${\rm ,} so
does the evaluation map $\overline D\ra A$.
\endproclaim

{\it Proof}. Consider the $C^*$-algebra $D[0,1]=D\otimes C[0,1]$ as a
$\Gamma-C([0,1]\times [0,1])$-algebra: the first copy of $C[0,1]$ in
$C([0,1]\times [0,1])=C[0,1]\otimes C[0,1]$ acting on $D$, the
second on $C[0,1]$. Let $U=\{(s,t)\in [0,1]\times [0,1]\,,
\;s\le t\,\}$ and $\overline D=D[0,1]_U$. Let
$C[0,1]$ act on $\overline D$, by $f\mt 1\otimes f$. In
this way, $\overline D$ is a $C[0,1]$ algebra. We have
$\overline D_{\{0\}}=A$. Moreover, let $$\varphi:[0,1]\times (0,1]\ra
\{(s,t)\in [0,1]\times (0,1]\,,\;s\le t\}$$ be a
homeomorphism of the form $(s,t)\mt (g(s,t),t)$ which is the
identity in a neighborhood of $\{0\}\times (0,1]\,$\footnote{\eg
$g(s,t)=\min (s,t(s+t)(1+t)^{-1})$.}. Extend $\varphi$ to a continuous map $[0,1]\times [0,1]\to [0,1]\times [0,1]$ by
the formula $\varphi(s,0)=(0,0)$.  Using
$\varphi$, we get an identification of the kernel of the evaluation
at $0$: $\overline D\ra A$ with $D(0,1]$; we therefore get
an asymptotic morphism $y:\quad 0\ra D(0,1]\ra \overline D\ra
A\ra 0$.

If $p_0$ admits an equivariant
completely positive lifting $h$, then $\pi\circ i\circ h$
is the lifting for the morphism $\overline D\ra A$, where
$i:D\ra D\otimes C[0,1]$ is the map $i:d\mt a\otimes 1$ and
$\pi :D\otimes C[0,1]\ra \overline D$ is the quotient map.

The algebra $\overline D$ is a $C(U)$ algebra. The asymptotic
morphism $p_0(y)$ is given by the exact sequence $0\ra
A(0,1]\ra \overline D_{\varphi (\{0\}\times [0,1])}\ra A\ra 0$; it is
the asymptotic morphism $0\ra A(0,1]\ra A[0,1]\ra A\ra 0$. The
asymptotic morphism $p_1(y)$ is given by the exact  sequence
$0\ra B(0,1]\ra \overline D_{\varphi (\{1\}\times [0,1])}\ra A\ra
0$. Now, $\varphi (1,s)=(s,s)$; the quotient map $D[0,1]\ra
\overline D_{\varphi (\{1\}\times [0,1])}$ composed with $i:D\ra
D[0,1]$ is  therefore a $C[0,1]$-linear isomorphism $D\ra
\overline D_{\varphi (\{1\}\times [0,1])}$, whence
$p_1(y)=x$.
\hfill\qed

\specialnumber{A.2}
\proclaim{Lemma}  Let $x:  \; 0\ra B(0,1]\ra D\ra
A\ra 0$  be an asymptotic morphism. Denote by $p_0:D\ra A$ and
$p_1:D\ra B$ the evaluation maps of the $C[0,1]$\/{\rm -}\/algebra
$D$.
\begin{itemize}
\ritem{a)} The morphism $p_0:D\ra A$ induces an isomorphism of
$E^\Gamma(A,D)$ onto
$E^\Gamma(A,A)$. The class of $x$ in $E^\Gamma(A,B)$ is
$p_{1*}((p_{0*})^{-1}(1^E_A))$.
\ritem{b)} If $p_0$ admits a completely positive equivariant lifting of norm $1$\/{\rm ,} 
the\break morphism
$p_0:D\ra A$ induces an isomorphism of $KK^\Gamma(A,D)$ onto
$KK^\Gamma(A,A)$. The class of $x$ in $KK^\Gamma(A,B)$ is
$p_{1*}((p_{0*})^{-1}(1^K_A))$.
\end{itemize}

\endproclaim

\demo{Proof}  The functor $B\lra E ^\Gamma(A,B)$ is `half exact' and the
functor $B\lra KK ^\Gamma(A,B)$ is `half exact' with respect to
exact sequences which admit a completely positive equivariant
lifting of norm $1$ ([BS]). Moreover, the ideal $B(0,1]$ is contractible in
an equivariant way; the first assertions in a) and b) follow. The second
ones are immediate consequences of Lemma A.1.
\enddemo

\specialnumber{A.3} \proclaim{Proposition} Let $A$ be a
$\Gamma$\/{\rm -}\/algebra. If the functor
$B\lra KK ^\Gamma(A,B)$ is {\rm `}\/half exact\/{\rm ',} then the homomorphism $\Phi:KK
^\Gamma(A,B)\lra E^\Gamma(A,B)$ is an isomorphism for any
$\Gamma$\/{\rm -}\/algebra $B$.
\endproclaim

\demo{Proof} We just have to show that the homomorphism $$\Phi:KK
^\Gamma(SA\otimes
\K,SB\otimes \K)\lra E^\Gamma(SA\otimes \K,SB\otimes \K)$$ is
an isomorphism. Since $$KK^\Gamma(SA\otimes \K,B)= KK^\Gamma(A,SB\otimes
\K),$$  the functor $B\lra KK^\Gamma(SA\otimes \K,B)$ is half exact. We
may therefore assume that there exist $\Gamma$-algebras
$A_1$ and $B_1$ such that $A=SA_1\otimes \K$ and
$B=SB_1\otimes \K$; therefore $E^\Gamma(A,B)=[[A,B]]$ and
$KK^\Gamma(A,B)=[[A,B]]_{\rm cp}$.

Let $x:\; 0\ra B(0,1]\ra D\ra A\ra 0$ be an asymptotic morphism
from $A$ to~$B$. Since the ideal $B(0,1]$ is contractible in
an equivariant way, the morphism $p_0:D\ra A$ induces an
isomorphism of $KK^\Gamma(A,D)$ onto
$KK^\Gamma(A,A)$. Let also $p_1:D\ra B$ be the evaluation map. To
the asymptotic morphism $x$, we associate the element $(p_1)_*\circ
(p_0)_*^{-1}(1_A^K)\in KK^\Gamma(A,B)$.

The $KK$-elements associated with homotopic asymptotic
morphisms coincide, since $KK^\Gamma$ is homotopy invariant. In this
way, we obtain a homomorphism
$\Psi:E^\Gamma(A,B)\ra KK^\Gamma(A,B)$.

Let $x:\; 0\ra B(0,1]\ra D\ra A\ra 0$ be an asymptotic morphism
from
$A$ to $B$ such that $p_0$ admits a completely positive
equivariant lifting of norm $1$. The class of $x$ in $KK^\Gamma(A,B)$
is $(p_1)_*\circ (p_0)_*^{-1}(1_A^K)$ (Lemma A.2.b). Therefore
$\Psi\circ \Phi$ is the identity (since $\Phi(1_A^E)=1_A^K$).

Let $x:\; 0\ra B(0,1]\ra D\ra A\ra 0$ be an asymptotic morphism
from $A$ to $B$. As $\Phi(1_A^E)=1_A^K$ and
$\Phi $ is natural in $B$, it follows that
$\Phi(\Psi(x))=\Phi((p_1)_*\circ
(p_0)_*^{-1}(1_A^K))=(p_1)_*\circ (p_0)_*^{-1}(\Phi(1_A^K))$ which
is the class of
$x$ in $E^\Gamma(A,B)$ (Lemma A.2.a). It follows that
$\Phi\circ \Psi$ is the identity.
\enddemo

Using Proposition 5.7 we obtain now:

\specialnumber{A.4} \proclaim{{C}orollary}\hskip-8pt Let $A$ be a nuclear proper
$\Gamma$\/{\rm -}\/algebra. Then $\Phi: KK^\Gamma(A,B)\break\ra E^\Gamma(A,B)$ is an
isomorphism.\endproclaim

{\it Remarks} A.5. 
a) Note that the converse of Proposition~A.3 is obviously true:
since the functor $B\lra E^\Gamma(A,B)$ is `half exact', if the
homomorphism $KK ^\Gamma(A,B)\lra E^\Gamma(A,B)$ is an
isomorphism for every $\Gamma$-algebra $B$, then the functor
$B\lra KK ^\Gamma(A,B)$ is `half exact'.
\vglue4pt
b) It follows that the Baum-Connes conjecture in $E$-theory of
[GHT] is equivalent to the one in $KK$-theory.

\AuthorRefNames [GHT ]

\end{document}

\centerline{(Received February 23, 2001)}
\bye